\documentclass[11pt]{amsart}
\usepackage{amsmath}
\usepackage[latin1]{inputenc}
\usepackage{amsfonts}
\usepackage{amssymb}

\usepackage{amsthm}
\usepackage{amscd}

\usepackage{pb-diagram}
\usepackage[mathscr]{eucal}
\usepackage[latin1]{inputenc}

\usepackage{amsmath}

\usepackage{amsfonts}

\usepackage{amssymb}

\usepackage{amsthm}

\usepackage{graphics}

\newcommand{\mk}{\medskip}

\newcommand{\ZZ}{\mathbb{Z}}

\newcommand{\CC}{\mathbb{C}}

\newcommand{\NN}{\mathbb{N}}

\newcommand{\Glie}{\mathfrak{g}}

\newcommand{\Yim}{\mathcal{Y}}

\newcommand{\demo}{\noindent {\it \small Proof:}\quad}

\newcommand{\U}{U}

\newcommand{\nc}{\newcommand}
\nc{\g}{{\mathfrak g}}
\nc{\ghat}{\widehat{\g}}
\nc{\mc}{\mathcal}
\nc{\ep}{\epsilon}
\nc{\la}{\lambda}
\nc{\Z}{{\mathbb Z}}
\nc{\C}{{\mathbb C}}
\nc{\on}{\operatorname}
\nc{\wt}{\widetilde}
\nc{\om}{\omega}
\nc{\ol}{\overline}
\nc{\wh}{\widehat}

\newtheorem{thm}{Theorem}[section]

\newtheorem{defi}[thm]{Definition}

\newtheorem{cor}[thm]{Corollary}

\newtheorem{prop}[thm]{Proposition}

\newtheorem{lem}[thm]{Lemma}

\newtheorem{conj}[thm]{Conjecture}

\newtheorem{rem}[thm]{Remark}

\usepackage[all]{xy}

\title{Langlands
  duality for finite-dimensional representations of quantum affine
  algebras}

\author[Edward Frenkel]{Edward Frenkel$^1$}\thanks{$^1$Supported 
in part by DARPA and AFOSR through the grant FA9550-07-1-0543 and by
Fondation Sciences mathématiques de Paris.}

\address{Department of Mathematics, University of California,
  Berkeley, CA 94720, USA}

\author[David Hernandez]{David Hernandez$^2$}\thanks{$^2$Supported
  partially by ANR through Project "G\'eom\'etrie et Structures
  Alg\'ebriques Quantiques"}

\address{D\'epartement de Math\'ematiques, Universit\'e Paris 7, 175 rue du Chevaleret, 75013 Paris, FRANCE}

\dedicatory{To V.B. Matveev on his 65th birthday}

\usepackage[all]{xy}
\begin{document}

\begin{abstract} 
We describe a correspondence (or duality) between the
$q$-characters of finite-dimensional representations of a quantum
affine algebra and its Langlands dual in the spirit of \cite{fr2, FH}.
We prove this duality for the Kirillov--Reshetikhin modules and their
irreducible tensor products.  In the
course of the proof we introduce and construct ``interpolating
$(q,t)$-characters'' depending on two parameters which interpolate
between the $q$-characters of a quantum affine algebra and its
Langlands dual.
\vskip 4.5mm



\end{abstract}

\maketitle

\tableofcontents

\pagestyle{myheadings}

\markboth{EDWARD FRENKEL AND DAVID HERNANDEZ}{LANGLANDS DUALITY FOR
  REPRESENTATIONS OF QUANTUM AFFINE ALGEBRAS}

\section{Introduction}

Let $\g$ be a simple Lie algebra and $\ghat$ the corresponding affine
Kac--Moody algebra. In \cite{fr2}, N. Reshetikhin and one of the
authors introduced a two-parameter deformed ${\mc W}$-algebra ${\mc
W}_{q,t}(\g)$. In the limit $t \to 1$ this deformed ${\mc W}$-algebra
becomes commutative and gives rise to the Grothendieck ring of
finite-dimensional representations of the quantum affine algebra
$\mathcal{U}_q(\ghat)$. (The precise relation between the two is
explained in \cite{fr2} and \cite{fr3}.) On the other hand, in the
limit when $q \to \ep$, where $\ep=1$ if $\g$ is simply-laced and $\ep
= \exp(\pi i/r)$, $r$ being the lacing number of $\g$, otherwise, this
algebra contains a large center. It was conjectured in \cite{fr2} that
it gives rise to the Grothendieck ring of $\U_q({}^L\ghat)$,
where $^L\ghat$ is the {\em Langlands dual Lie algebra} of $\ghat$. By
definition, the Cartan matrix of $^L\ghat$ is the transpose of the
Cartan matrix of $\ghat$, so that $^L\ghat$ is a twisted affine
algebra if $\g$ is non-simply laced.

Thus, it appears that ${\mc W}_{q,t}(\g)$ {\em interpolates} between
the Grothendieck rings of finite-dimensional representations of
quantum affine algebras associated to $\ghat$ and $^L\ghat$. This
suggests that these representations should be related in some
way. Examples of such a relation were given in \cite{fr2}, but general
understanding of this phenomenon has been lacking. The goal of this
paper is to elucidate and provide further evidence for this duality.

The finite-dimensional analogue of this duality has been studied in
our previous paper \cite{FH}, in which we have conjectured (and
partially proved) the existence of a correspondence, or duality,
between finite-dimensional representations of the quantum groups
$\U_q(\g)$ and $\U_q({}^L \g)$.\footnote{We have learned from
  K. McGerty that in the meantime he has been able to prove one of the
  conjectures of \cite{FH}, see \cite{M}. After this paper was
  published, we learned from C. Lecouvey that the statement of
  this conjecture was proved earlier by P. Littelmann in \cite{pl}.}
This duality may in fact be extended uniformly to integrable
representations of quantized enveloping algebras associated to
Kac--Moody algebras. But quantized enveloping algebras associated to
the affine Kac--Moody algebras (quantum affine algebras for short)
have another important class of representations; namely, the
finite-dimensional representations.  In this paper we describe a
Langlands type duality for these representations.

In this context the Langlands duality was first observed in \cite{fr2,
  fr3} using the so-called ``$q$-characters'' of finite-dimensional
representations of quantum affine algebras. The theory of
$q$-characters has been developed for untwisted quantum algebras in
\cite{fr3} and for twisted quantum affine algebras (which naturally
appear in the Langlands dual situation) in \cite{hgen}.  The
$q$-character of a simple representation characterizes its isomorphism
class.

In this paper we conjecture a precise relation between the
$q$-characters of finite-dimensional representations of dual quantum
affine algebras $\U_q(\ghat)$ and $\U_q({}^L\ghat)$. Namely, we
conjecture that for any finite-dimensional representation $V$ of
$\U_q(\wh{\Glie})$ there exists an {\it interpolating
$(q,t)$-character}, a polynomial which interpolates between the
$q$-character of $V$ and the $t$-character of a certain representation
of the Langlands dual algebra $\U_t({}^L\wh{\Glie})$, which we call
{\em dual to} $V$ (we discuss in which sense it is unique). Moreover, we prove this
conjecture for an important class of representations, the {\em
Kirillov--Reshetikhin modules} and their irreducible tensor products. 

The $q$-characters are important in the study of integrable models of
statistical mechanics (see, e.g., \cite{fr2,fr3}), and therefore their duality
indicates the existence of duality between the models associated to
two Langlands dual affine Lie algebras. The existence of interpolating
$(q,t)$-characters is closely related to \cite[Conjecture 1]{fr2},
which also states the existence of interpolating expressions, but of a
different kind. They are elements of a two-parameter non-commutative
algebra (in fact, a Heisenberg algebra), whereas the interpolating
characters that we introduce here are elements of a commutative
algebra. It would be interesting to understand a precise relation
between the two pictures.

We refer the reader to the Introduction of \cite{FH} for a discussion
of a possible link between our results on the duality of
finite-dimensional representations of $\U_q(\ghat)$ and
$\U_t({}^L\ghat)$ and the geometric Langlands correspondence. This
link was one of the motivations for the present paper.

Let us note that the technique and methods in the present paper are
not generalizations of \cite{FH}, but are new as we use the
``rigidity'' provided by the appearance of the spectral parameters in
the context of quantum affine algebras. This allows us to construct
the interpolating $(q,t)$-characters (which have {\em a priori} no
clear analogues in finite types). Another difference is that instead
of a projection from a weight lattice to the dual weight lattice, we
introduce interpolating maps $\alpha(q,t), \beta(q,t)$ in the
characters.  These maps ``kill'' some of the terms when we specialize
to the Langlands dual situation. Thus we obtain a much finer form of
duality in the affine case than in the finite-dimensional case.

\medskip

The paper is organized as follows: in Section \ref{usual} we recall
the Langlands duality for quantum groups of finite type from
\cite{FH}. Then we state consequences of the results of the present
paper in terms of the ordinary characters (Theorem \ref{thusual}). In
Section \ref{double} we give a general conjecture about the duality at
the level of $q$-characters. We state and start proving the main
result of the present paper (Theorem \ref{thdual}) in the double-laced
cases; namely, that the Kirillov--Reshetikhin modules and their irreducible
tensor products satisfy the
Langlands duality. The end of the proof uses results of Section
\ref{interpolate} where interpolating $(q,t)$-characters are
constructed in a systematic way (Theorem \ref{con}). The triple-laced
is treated in Section \ref{typeg} (Theorem \ref{virtg} and Theorem
\ref{dualg}) to complete the picture. In Section \ref{reverse}, we
describe a reverse Langlands duality from twisted quantum affine
algebras to untwisted quantum algebras and we prove analogous results
for this duality (Theorem \ref{virtrev} and Theorem \ref{dualrev}).

\section{Duality for the ordinary characters}\label{usual}

Although most of the results of the present paper involve
$q$-characters, some consequences of our results may be stated purely
in terms of the ordinary characters. We explain these results in this
Section as well as some motivations and results from \cite{FH}.

Let $\g$ be a finite-dimensional simple Lie algebra and $\U_q(\g)$ the
corresponding quantum group (see, e.g., \cite{cp}). We denote $r =
\text{max}_{i\in I}(r_i)$, where $I$ is the set of vertices of the
Dynkin diagram of $\g$ and the $r_i$ are the corresponding
labels. This is the lacing number of $\g$ (note that it was denoted by
$r^\vee$ in \cite{fr2,fr3}).  In some particular cases, we will not
make the choice $\text{min}_{i\in I}(r_i) = 1$ (that is we multiply
the standard labels by a coefficient).

The Cartan matrix of $\g$ will be denoted by $C=(C_{i,j})_{i,j \in
I}$. By definition, the Langlands dual Lie algebra $^L\g$ has the
Cartan matrix $C^t$, the transpose of the Cartan matrix $C$ of $\g$.

Let $$P=\sum_{i \in I} \ZZ\om_i$$ be the weight lattice of $\Glie$ and
$P^+ \subset P$ the set of dominant weights. For $i\in I$ let
$r_i^\vee = 1 + r - r_i$ and consider the sublattice
\begin{equation}    \label{Pprime}
P' = \sum_{i\in I} r_i^\vee \ZZ\om_i \subset P.
\end{equation}
Let $$P^L=\sum_{i \in I} \ZZ\check\om_i$$ be the weight lattice of
$^L\g$. Consider the map $\Pi : P\rightarrow P^L$ defined by
$$\Pi(\lambda) = \sum_{i\in I}\lambda(\check\alpha_i)(r_i^\vee)^{-1}
\check\omega_i$$ if $\lambda \in P'$ and $\Pi(\lambda) = 0$,
otherwise. Clearly, $\Pi$ is surjective.

Let $\on{Rep} \g$ be the Grothendieck ring of finite-dimensional
representations of $\g$. We have the character homomorphism
$$
\chi: \on{Rep} \g \to \Z[P] = \Z[y_i^{\pm 1}],
$$
where $y_i = e^{\om_i}$. It sends an irreducible representation
$L(\la)$ of $\g$ with highest weight $\la \in P^+$ to its character,
which we will denote by $\chi(\la)$.  We denote the character
homomorphism for $^L\g$ by $\chi^L$. We use the obvious partial
ordering $\preceq$ on polynomials. It was proved in \cite{FH} that
for any $\lambda \in P^+$, $\Pi(\chi(\lambda))$ is in the image of
$\chi^L$. Moreover, we have the following:

\begin{thm}\label{pos}\cite{FH} For any $\lambda \in P^+$,
  $\Pi(\chi(\lambda))\succeq \chi^L(\Pi(\lambda))$.
\end{thm}

Let $q,t\in\CC^\times$ be such that $q^\ZZ\cap t^\ZZ = \{1\}$. Let
$\U_q(\wh{\Glie})$ be an untwisted quantum affine algebra which is not
Langlands self-dual.  Let $V$ be a simple finite-dimensional
representation of $\U_q(\wh{\Glie})$ which is of highest weight
$\lambda$ in $P'$ when viewed as a $\U_q(\Glie)$-module. We conjecture
the following.

\begin{conj}\label{usualchar} There exists an irreducible representation
  $V^L$ of $\U_t(\wh{\Glie}^L)$ of highest weight $\Pi(\lambda)$ such
  that $\Pi(\chi(V))\succeq \chi^L(V^L)$.
\end{conj}

Note that the Langlands dual representation $V^L$ is not necessarily
unique. Uniqueness statements will be discussed later in a more precise
form of Langlands duality.

As a consequence of the results of the present paper, we will prove the following.

\begin{thm}\label{thusual} The statement of Conjecture \ref{usualchar}
  is satisfied for any Kirillov--Reshetikhin module $V$ over
  $\U_q(\wh{\Glie})$, with the Langlands dual representation $V^L$ 
  a Kirillov--Reshetikhin module over $\U_t({}^L\wh{\Glie})$.
\end{thm}

Note also that in contrast to \cite{FH}, we use $t$ and not $-t$ for the
quantization parameter of the Langlands dual quantum algebra. This is
just a consequence of a different choice of normalization made in the
present paper.

The following conjecture of \cite{FH} has been proved by K. McGerty in
\cite{M}\footnote{After this paper was
  published, we learned from C. Lecouvey that the statement of
  this conjecture was proved earlier by P. Littelmann in \cite{pl}.}: for any $\lambda \in P^+$, $\Pi(\chi(\lambda))$ is the
character of an actual (not only virtual) representation of
${}^L\Glie$. Therefore it is natural to make the following.

\begin{conj}
$\Pi(\chi(L))$ is the character of a representation of
$\U_t({}^L\wh{\Glie})$.
\end{conj}

Again, this representation of $\U_t({}^L\wh{\Glie})$ is
not unique, but it is unique as a $\U_t({}^L\Glie)$-module. However it
is not necessarily simple as a $\U_t({}^L\Glie)$-module. As an
example, for a $5$-dimensional fundamental representation of
$\U_q(C_2^{(1)})$, the Langlands dual representation decomposes into a
sum of a $4$-dimensional fundamental representation and the trivial
representation of $\U_t(A_3^{(2)})$ (see the corresponding
$q$-characters in Section \ref{exa}).

\section{Double-laced cases}\label{double}

In this section we suppose that the lacing number $r$ is equal to $2$
(the case $r = 3$ will be treated in Section \ref{typeg}).  We will
exclude from consideration the Langlands self-dual quantum affine Lie
algebras (affinizations of simply-laced ones and those
$A_{2n}^{(2)}$).

We have $I = I_1\sqcup I_2$ where $I_k = \{i\in I| r_i = k\}$. For
$i,j\in I$, we denote $i\sim j$ if $C_{i,j} < 0$. We can choose $\phi
: I\rightarrow \{1 , 0\}$ such that $i\sim j \Rightarrow \phi(i) +
\phi(j) = 1$ and $C_{i,j} = -2\Rightarrow \phi(i) = 1$.

\subsection{Reminder on $q$-characters and their twisted analogues} We
recall the notion of $q$-characters first introduced in \cite{fr3} for
untwisted quantum affine algebras (see \cite{ch} for a recent survey)
and generalized in \cite{hgen} to the twisted cases.

The $q$-character homomorphism \cite{fr3} is an injective ring morphism
$$\chi_q : \text{Rep}(\U_q(\wh{\Glie}))\rightarrow \Yim_q =
\ZZ[Y_{i,a}^{\pm 1}]_{i\in I, a\in q^\ZZ}$$ (without loss of
generality, we restrict ourselves to the tensor subcategory of
finite-dimensional representations whose $q$-characters are in
$\Yim_q$). By removing the spectral parameter $a$, that is to say by
replacing each $Y_{i,a}$ by $y_i$, we recover the usual character map
for the $\U_q(\Glie)$-module obtained by restriction of
$\U_q(\wh{\Glie})$-module. In particular, each monomial has a weight
which is an element of $P$. For $i\in I$, let $q_i = q^{r_i}$.

\begin{thm}\cite{fm} We have 
$$\on{Im}(\chi_q) = \bigcap_{i\in I}\mathfrak{K}_{i,q},$$ where
$\mathfrak{K}_{i,q} = \ZZ[Y_{j,a}^{\pm 1}, Y_{i,a}(1 +
A_{i,aq_i}^{-1})]_{j\neq i,a\in q^\ZZ}$ and
$$A_{i,a} = Y_{i,aq_i^{-1}}Y_{i,aq_i} \times \prod_{j\in I,C_{j,i} =
    - 1}Y_{j,a}^{-1} \times \prod_{j\in I,C_{j,i} = -
    2}Y_{j,aq^{-1}}^{-1} Y_{j,aq}^{-1}.$$
\end{thm}

A monomial in $\Yim_q$ is called {\em dominant} if it is a product of
positive powers of the $Y_{i,a}$ (for $i\in I, a\in q^{\ZZ}$). A
simple $\U_q(\wh{\Glie})$-module is characterized by the highest
monomial (in the sense of its weight) of its $q$-character (this is
equivalent to the data of the Drinfeld polynomials, see
\cite{cp}). This monomial is dominant. Any element of
$\text{Im}(\chi_q)$ is characterized by the list of its dominant
monomials. A $\U_q(\wh{\Glie})$-module is said to be {\em
affine-minuscule} if its $q$-character has a unique dominant monomial.

\begin{defi} A Kirillov--Reshetikhin (KR) module of $\U_q(\wh{\Glie})$ is a
  simple module with the highest monomial of the form
  $Y_{i,a}Y_{i,aq_i^2}\cdots Y_{i,aq_i^{2(k-1)}}$.
\end{defi}

We have the following result which is due to H. Nakajima \cite{Nab,
Nad} in the simply-laced case and \cite{hcr} in general (note that for
$k=1$ this was proved in \cite{fm} in the untwisted case).

\begin{thm}\label{krun} The KR modules of
  $\U_q(\wh{\Glie})$ are affine-minuscule.
\end{thm}

Now let us look at the Langlands dual situation, i.e.,
finite-dimensional representations of the twisted quantum affine
algebra $\U_t({}^L\wh{\Glie})$.  We set $\epsilon = e^{i\pi/2}$
and $I_1^\vee = I_2, I_2^\vee = I_1$.

The twisted $t$-character morphism \cite{hgen} is an injective ring
homomorphism (we work in a subcategory defined as in the untwisted
case)
$$\chi_t^L : \text{Rep}({}^L\U_t(\wh{\Glie}))\rightarrow \Yim_t^L =
\ZZ[Z_{i,a^{r_i^\vee}}^{\pm 1}]_{a\in \epsilon^\ZZ t^\ZZ, i\in I}.$$

\begin{thm}\cite{hgen} We have 
$$\on{Im}(\chi_t^L) = \bigcap_{i\in I}\mathfrak{K}_{i,t}^L,$$ 
where 
$$\mathfrak{K}_{i,t}^L = \ZZ[Z_{j,a^{r_j^\vee}}^{\pm 1}, Z_{i,a^{r_i^\vee}}(1 +
  B_{i,(at)^{r_i^\vee}}^{-1})]_{j\neq i,a\in \epsilon^\ZZ t^\ZZ}$$
and
$$B_{i,a} = \begin{cases} Z_{i,at^2}Z_{i,at^{-2}}\times \prod_{j\sim
i|j \in I_2^\vee}Z_{j,a}^{-1}\times \prod_{j\sim i| j \in I_1^\vee}\prod_{a'\in
\epsilon^\ZZ t^\ZZ|(a')^2 = a}Z_{j,a'}^{-1} & \text{ if $i \in I_2^\vee$,}
\\Z_{i,at}Z_{i,at^{-1}}\times \prod_{j\sim i|j \in
  I_2^\vee}Z_{j,a^2}^{-1}\times \prod_{j\sim i| j\in
  I_1^\vee}Z_{j,a}^{-1}&\text{ if $i \in I_1^\vee$.}\end{cases}$$
\end{thm}

Note that a special definition should be used \cite{hgen} for the
$B_{i,a}$ in the case of type $A_{2n}^{(2)}$, but this case is not
considered here.

We have the notions of dominant monomial, affine-minuscule module and
KR module as in the untwisted case. Any element of
$\text{Im}(\chi_t^L)$ is again characterized by its dominant
monomials and we have

\begin{thm}\label{krdeux}\cite{hgen} The KR modules
  of $\U_t({}^L\wh{\Glie})$ are affine-minuscule.
\end{thm}

\subsection{The interpolating $(q,t)$-character ring}

We first treat study the duality from untwisted quantum affine
algebras to twisted quantum affine algebras. The reverse Langlands
duality will be treated later.

We introduce the {\em interpolating $(q,t)$-characters}, which
interpolate between $q$-characters of an untwisted quantum affine
algebra and the twisted $t$-characters of its Langlands dual. To do
it, we first need to define an interpolating ring for the target rings
of $q$-and $t$-character homomorphisms.

We also need the function $\alpha(q,t)$ such that $\alpha(q,1) =
1$ and $\alpha(\epsilon, t) = 0$ defined in
\cite{fr2,fr3} (see also \cite{FH} for an elementary natural way to
introduce it in the framework of current algebras) by the formula
$$\alpha(q,t) = \frac{(q + q^{-1})(qt - q^{-1}t^{-1})}{q^2t -
q^{-2}t^{-1}}.$$

Let $\mathcal{C} = q^\ZZ t^\ZZ$. Consider the ring 
$$\Yim_{q,t} = \ZZ[W_{i,a}^{\pm 1}, \alpha Y_{i,a}^{\pm 1},
  \alpha]_{i\in I, a\in\mathcal{C}}\subset \ZZ[Y_{i,a}^{\pm
    1},\alpha]_{i\in I, a\in\mathcal{C}},$$
\begin{equation*}
\begin{split}
\text{ where }W_{i,a} = \begin{cases} Y_{i,a}&\text{ if $i\in I_2$,}
                     \\Y_{i,aq^{-1}}Y_{i,aq}&\text{ if $i\in I_1$.}
\end{cases}
\end{split}
\end{equation*}
For $a\in \mathcal{C}$, we will use the following identification for
$i\in I_1$ and $j\in I_2$:
\begin{equation}\label{ident} Y_{i,a}Y_{i,-a} = Z_{i, a^2
    (-1)^{\phi(i)}}\text{   and   }Y_{j,a} =
  Z_{j,a(-1)^{\phi(i)}}.
\end{equation}

We then have surjective specialization maps, respectively, at $t = 1$
and $q = \epsilon$,
$$\Pi_q : \Yim_{q,t}\rightarrow \ZZ[Y_{i,a}^{\pm 1}]_{i\in I, a\in
  q^\ZZ} = \Yim_q,$$ $$\Pi_t : \Yim_{q,t}\rightarrow
  \ZZ[Z_{i,a^{r_i^{\vee}}}^{\pm 1}]_{i\in I, a\in \epsilon^\ZZ t^\ZZ}
  = \Yim_t^L.$$

We have the ideals
$$\text{Ker}(\Pi_q) = \langle (\alpha - 1), (W_{i,a} - W_{i,at}),
\alpha(Y_{i,a} - Y_{i,at}) \rangle_{i\in I,a\in\mathcal{C}},$$
$$\text{Ker}(\Pi_t) = \langle \alpha, (W_{i,aq} - W_{i,a\epsilon})
  \rangle_{i\in I, a\in\mathcal{C}}.$$ 
  Then we have the following:

\begin{lem} The ideal $\text{Ker}(\Pi_q)\cap \text{Ker}(\Pi_t)$ is
  generated by the elements 
$$\alpha (\alpha - 1)\text{ , }\alpha (Y_{i,a} -
  Y_{i,at})\text{ , }(\alpha - 1) (W_{i,aq} - W_{i,a\epsilon})\text{ ,
  }(W_{i,a} - W_{i,at})(W_{j,bq} - W_{j,b\epsilon}),$$
for $i,j\in I$, $a,b\in\mathcal{C}$.
\end{lem}

\demo First, the ideal $\mathfrak{I}$ generated by these elements is clearly
included the intersection $\text{Ker}(\Pi_q)\cap \text{Ker}(\Pi_t)$
and so we work modulo $\mathfrak{I}$. We denote by $\equiv$ the equality
modulo $\mathfrak{I}$. Now consider an element $\chi$ in the
intersection. It is of the form
$$\chi \equiv (\alpha - 1)\chi(q,t) + \sum_{i\in I, l,r\in\ZZ}
(W_{i,q^lt^r} - W_{i,q^l t^{r + 1}})\chi_{i,l,r}(q,t)$$
$$\equiv (\alpha - 1) \chi(\epsilon, t)+ \sum_{i,l,r}  (W_{i,q^l t^r} -
W_{i,q^l t^{r + 1}}) \chi_{i,l,r}(\epsilon,t).$$
If we evaluate at $q = \epsilon$, we get
$$\chi(\epsilon,t) \equiv \sum_{i,l,r}  (W_{i,\epsilon^l t^r} -
W_{i,\epsilon^l t^{r + 1}}) \chi_{i,l,r}(\epsilon,t),$$
and so
$$\chi\equiv \sum_{i,l,r}  (W_{i,q^l t^r} - W_{i,q^l t^{r + 1}} + (\alpha -
1)(W_{i,\epsilon^l t^r} - W_{i,\epsilon^l t^{r + 1}}))
\chi_{i,l,r}(\epsilon,t)$$
As $(\alpha - 1)(W_{i,\epsilon^l t^r} - W_{i,\epsilon^l t^{r + 1}}) \equiv
-(W_{i,q^l t^r} - W_{i,q^l t^{r + 1}})$, we get
$$\sum_{i,l,r}  (W_{i,q^l t^r} - W_{i,q^l t^{r + 1}} - (W_{i,q^l t^r}
- W_{i,q^l t^{r + 1}})) \chi_{i,l,r}(\epsilon,t) \equiv 0.$$
This concludes the proof.
\qed

We will work in the ring
$$\wt{\Yim}_{q,t} = \Yim_{q,t}/(\text{Ker}(\Pi_q)\cap
\text{Ker}(\Pi_t)).$$ Note that $\wt{\Yim}_{q,t}$ has zero divisors
as $\alpha^2 = \alpha$ in $\wt{\Yim}_{q,t}$. 

By a monomial in $\wt{\Yim}_{q,t}$ we will understand an element $m$
of the form $(\lambda + \mu \alpha)M$, where $\lambda,\mu\in\ZZ$ and
$M$ is a monomial in the $Y_{j,a}^{\pm 1}$.  Note that a monomial may
be written in various way as for example $\alpha Y_{i,a} = \alpha
Y_{i,at}$ and $(1 - \alpha)Y_{i,aq^4} = (1 - \alpha)Y_{i,a}$.  A
monomial is said to be $i$-dominant if it can be written by using only
the $\alpha$, $Y_{i,a}$ and $Y_{j,a}^{\pm 1}$ where $j\neq i$. Let
$B_i$ be the set of $i$-dominant monomials and for $J\subset I$, let 
$B_J = \cap_{j\in J} B_j$. Finally, $B = B_I$ is the set of dominant
monomials.

\subsection{Subalgebras of $\wt{\Yim}_{q,t}$}

\begin{defi}For $i\in I$ and $a\in \mathcal{C}$ we define
$$\wt{A}_{i,a} = Y_{i,a(q_it)^{-1}}Y_{i,aq_it} \times \prod_{j\in
I,C_{j,i} = - 1}Y_{j,a}^{-1}\times \prod_{j\in I,C_{j,i} = -
2}Y_{j,aq^{-1}}^{-1} Y_{j,aq}^{-1}.$$
\end{defi}
Note that the definition of $\wt{A}_{i,a}$ is not symmetric in
$q,t$. For $i\in I_2, a\in\mathcal{C}$ we have $\wt{A}_{i,a}^{\pm
1}\in \Yim_{q,t}$, and for $i\in I_1$ we have $\alpha
\wt{A}_{i,a}^{\pm 1}\in\Yim_{q,t}$ and
$(\wt{A}_{i,aq^{-1}}\wt{A}_{i,aq})^{\pm 1}\in\Yim_{q,t}$. But the
specialization maps $\Pi_q, \Pi_t$ can be applied to any
$\wt{A}_{i,a}$ and we have the following:

\begin{lem}\label{identa} We have $\Pi_q(\wt{A}_{i,a}) =
  A_{i,\Pi_q(a)}$ for $i\in I, a\in\mathcal{C}$.

We have $\Pi_t((\wt{A}_{i,aq^{-1}}\wt{A}_{i,aq})) =
B_{i,(\Pi_t(a))^2(-1)^{\phi(i)}}$ for $i\in I_1, a\in\mathcal{C}$.

We have $\Pi_t(\wt{A}_{i,a}) =
B_{i,-\Pi_t(a)(-1)^{\phi(i)}}$ for $i\in I_2, a\in\mathcal{C}$.
\end{lem}

\demo The first point is clear.

Let $a' = \Pi_t(a)$. For $i\in I_1$, the specialization of
$\wt{A}_{i,aq^{-1}}\wt{A}_{i,aq}$ at $q = \epsilon$ is
\begin{equation*}
\begin{split}
&(Y_{i,-a' t^{-1}} Y_{i,a' t^{-1}})(Y_{i,a't}Y_{i,-a' t}) \times
\prod_{j\in I_1,j\sim i}(Y_{j,a'\epsilon}Y_{j,-a'\epsilon})^{-1}
\times \prod_{j\in I_2, j\sim
i}Y_{j,a'\epsilon}^{-1}Y_{j,-a'\epsilon}^{-1} \\=
&Z_{i,(a')^2t^{-2}(-1)^{\phi(i)}}Z_{i,(a')^2t^2(-1)^{\phi(i)}} \times
\prod_{j\in I_1, j\sim i}Z_{j,(a')^2(-1)^{\phi(i)}}^{-1}
\times \prod_{j\in I_,  j\sim i} Z_{j,a'\epsilon}^{-1}Z_{j,-a'\epsilon}^{-1}.
\end{split}
\end{equation*}
Note that if there is $j\in I_2$ such that $j\sim i$, by definition of
$\phi$ we have $\phi(j) = 0$ and $\phi(i) =1$. That is why there is no
$\phi$ in the last factor of the product.

For $i\in I_2$, the specialization of $\wt{A}_{i,a}$ at $q = \epsilon$
is
$$Y_{i,-a't^{-1}}Y_{i,-a't} \times \prod_{j\in I_2,j\sim
i}Y_{j,a'}^{-1}\times \prod_{j\in I_1, j\sim
i}(Y_{j,-a'\epsilon}Y_{j,a'\epsilon})^{-1}$$
$$= Z_{i,-a't^{-1}(-1)^{\phi(i)}}Z_{i,-a't(-1)^{\phi(i)}}\times
\prod_{j\in I_2,j\sim i}Z_{j,-a'(-1)^{\phi(i)}}^{-1}\times \prod_{j\in
I_1, j\sim i} Z_{j,(a')^2}^{-1} .$$ Note that if there is $j\in I_1$
such that $j\sim i$, by definition of $\phi$ we have $\phi(j) = 1$ and
$\phi(i) = 0$. That is why there is no $\phi$ in the last factor of
the product.  \qed

For $i\in I_1$ consider the subalgebra of $\Yim_{q,t}$
$$\mathfrak{K}_{i, q, t} = \ZZ[W_{i,a}(1 + \alpha
\wt{A}_{i,aq^2t}^{-1} + \wt{A}_{i,aq^2t}^{-1}\wt{A}_{i,at}^{-1}),
\alpha Y_{i,a}(1 + \wt{A}_{i,aqt}^{-1}), W_{j,a}^{\pm 1}, \alpha
Y_{j,a}^{\pm 1}, \alpha]_{a\in\mathcal{C}, j\neq i},$$ and for $i\in
I_2$
$$\mathfrak{K}_{i,q, t} = \ZZ[Y_{i,a}(1 + \wt{A}_{i,aq^2t}^{-1}),
  W_{j,a}^{\pm 1}, \alpha Y_{j,a}^{\pm 1}, \alpha]_{a\in \mathcal{C},
  j\neq i}.$$ 
Then we have the following:

\begin{lem}\label{special} For $i\in I$, we have 
$\Pi_q(\mathfrak{K}_{i,q,t})= \mathfrak{K}_{i,q}$ and
  $\Pi_t(\mathfrak{K}_{i,q,t}) = \mathfrak{K}_{i,t}^L$.
\end{lem}

\demo For $i\in I_1$, we have
$$\Pi_q(\mathfrak{K}_{i,q,t}) = \ZZ[Y_{i,a}(1 + A_{i,aq}^{-1}),
Y_{j,a}^{\pm 1}]_{a\in q^\ZZ, j\neq i} = \mathfrak{K}_{i,q},$$
$$\Pi_t(\mathfrak{K}_{i,q,t}) = \ZZ[Y_{i,a}Y_{i,-a}(1 +
(A_{i,-\epsilon at}A_{i,\epsilon at})^{-1}), (Y_{j,a}Y_{j,-a})^{\pm
1}, Y_{k,a}^{\pm 1}]_{a\in \epsilon^\ZZ t^\ZZ, j\in I_1-\{i\}, k\in
I_2} = \mathfrak{K}_{i,t}^L,$$ as by Lemma \ref{identa} we have
$$Y_{i,a}Y_{i,-a}(1 + (A_{i,-\epsilon at}A_{i,\epsilon at})^{-1}) =
Z_{i,a^2(-1)^{\phi(i)}}(1 + B_{i,a^2t^2(-1)^{\phi(i)}}^{-1}).$$ For
$i\in I_2$, we have
$$\Pi_q(\mathfrak{K}_{i,q,t}) = \ZZ[Y_{i,a}(1 + A_{i,aq^2}^{-1}),
Y_{j,a}^{\pm 1}]_{a\in q^\ZZ, j\neq i} = \mathfrak{K}_{i,q},$$
$$\Pi_t(\mathfrak{K}_{i,q,t}) = \ZZ[Y_{i,a}(1 + A_{i,-at}^{-1}),
(Y_{j,a}Y_{j,-a})^{\pm 1}, Y_{k,a}^{\pm 1}]_{a\in \epsilon^\ZZ t^\ZZ,
j\in I_1, k\in I_2-\{i\}} = \mathfrak{K}_{i,t}^L,$$ as by Lemma
\ref{identa} we have
$$Y_{i,a}(1 + A_{i,-at}^{-1}) = Z_{i,a (-1)^{\phi(i)}} (1 + B_{i, a
(-1)^{\phi(i)} t}^{-1}).$$ \qed

We use the same notation $\mathfrak{K}_{i,q,t}$ for the image of the
subalgebras $\mathfrak{K}_{i,q,t}$ in $\wt{\Yim}_{q,t}$.
For $J\subset I$ we define $\mathfrak{K}_J =
\bigcap_{j\in J} \mathfrak{K}_j$ and we set $\mathfrak{K} =
\mathfrak{K}_I\subset \tilde{\Yim}_{q,t}$.

\subsection{Main conjecture and main theorem}

Let us define an analogue of $P'\subset P$,
\begin{equation*}
\Yim_q' = \ZZ[Y_{i,a}^{\pm 1}]_{i \in I_2,a\in q^\ZZ}\otimes
\ZZ[(Y_{i,aq}Y_{i,aq^{-1}})^{\pm 1}]_{i \in I_1,a\in q^\ZZ}\subset
\Yim_q.
\end{equation*}

We consider $\wh{\Pi} : \Yim_q \rightarrow \Yim_q'$ the projection on
$\Yim_q'$ whose kernel is generated by monomials not in $\Yim_q'$.

Let $M\in\Yim_q'$ be a dominant monomial and 
$V$ the corresponding irreducible representation of $\U_q(\wh{\Glie})$. 
A representation $V^L$ of $\U_t({}^L\wh{\Glie})$ is said to be
{\em Langlands dual to} $L(M)$ if there is 
a dominant monomial $\wt{M}\in\wt{\Yim}_{q,t}\setminus\alpha \wt{\Yim}_{q,t}$
and $\wt{\chi}_q\in \mathfrak{K}\cap \wt{M}\ZZ[\wt{A}_{i,a}^{-1},
  \alpha]_{i\in I, a\in\mathcal{C}}$ such that
$$\Pi_t(\wt{\chi}_q) = \chi_q^L(V^L)\text{ and }\Pi_q(\wt{\chi}_q) = \chi_q(V).$$
Besides, we say that $\wt{\chi}_q$ is an {\em interpolating $(q,t)$-character} of $V$.

A given representation $V$ may have different Langlands dual
representations (for example, obtained by a shift of the spectral
parameter by $t^N$, that is by replacing each $Z_{i,a^{r_i^\vee}}$ by
$Z_{i,(at^N)^{r_i^\vee}}$ in $\chi_t^L(V^L)$). Following \cite{cp2,
  hl}, we will call a representation which cannot be factorized as a
tensor product of non-trivial representations, a {\em prime representation}.

\begin{conj}\label{char} Any irreducible representation $V$ of $\U_q(\wh{\Glie})$ has a
  Langlands dual representation $V^L$. Moreover, if $V$ is prime,
  $V^L$ is unique up to a shift of spectral parameter.
\end{conj}

If $V$ is not prime, the uniqueness statement does not necessarily hold 
(see Remark \ref{uni}).
Conjecture \ref{char} implies Conjecture \ref{usualchar} as the
condition $\wt{M}\in\wt{\Yim}_{q,t}\setminus \alpha \wt{\Yim}_{q,t}$
implies that the highest weight of $V^L$ is given
by the highest weight of $V$. The following is the main result of this paper.

\begin{thm}\label{thdual} Let $V$ be a KR module over $\U_q(\wh{\Glie})$ or an irreducible
  tensor product of KR modules. Then $V$ has a Langlands dual
  representation. Moreover, the Langlands dual representation of a KR
  module over $\U_q(\wh{\Glie})$ is a KR module over $\U_t({}^L\wh{\Glie})$.
\end{thm}

To prove this Theorem, we will use the affine-minuscule property of the KR modules. 

\subsection{Examples}\label{exa}

Let us give some examples of interpolating $(q,t)$-characters which will
be useful in the following proofs.

First consider the type $A_1$ with $r = 1$. We choose $\phi(1) = 0$
and we have the following.
$$\xymatrix{
   Y_{1,1}Y_{1,q^2} \ar[d]^{1,q^3t}             &Y_{1,1}Y_{1,q^2} \ar[d]^{1,q^3}     & Z_{1,1}\ar[dd]^{1,t^2}
\\ \alpha Y_{1,1}Y_{1,q^4t^2}^{-1}\ar[d]^{1,qt} &Y_{1,1}Y_{1,q^4}^{-1}\ar[d]^{1,q}   &
\\ Y_{1,q^2t^2}^{-1}Y_{1,q^4t^2}^{-1}             &Y_{1,q^2}^{-1}Y_{1,q^4}^{-1}      & Z_{1,t^4}^{-1}            
}$$
Here we use diagrammatic formulas for (interpolating) $q$-characters
as defined in \cite{fr3}. The left term in the interpolating
$q$-character, and then we have the respective specializations at $t =
1$ and $q = \epsilon$.

Consider the type $A_1$ with $r = 2$. We choose $\phi(1) = 0$ and we have the following.
$$\xymatrix{
   Y_{1,1} \ar[d]^{1,q^2t}    &Y_{1,1}       \ar[d]^{1,q^2}   &Z_{1,1}\ar[d]^{1,t}
\\ Y_{1,q^4t^2}^{-1}          &Y_{1,q^4}^{-1}                 &Z_{1,t^2}^{-1}
}$$
Next, consider the type $A_2$ with $r = 1$. We choose $\phi(1) = 0$, $\phi(2) = 1$ and we have the following.
$$\xymatrix{
   Y_{1,1}Y_{1,q^2} \ar[d]^{1,q^3t}        &&Z_{1,1}\ar[dd]^{1,t^2}
\\ \alpha Y_{1,1}Y_{1,q^4t^2}^{-1}Y_{2,q^3t}\ar[d]^{1,qt}\ar[dr]^{2,q^4t^2} &&
\\ Y_{1,q^2t^2}^{-1}Y_{1,q^4t^2}^{-1}Y_{2,q^3t}Y_{2,qt}\ar[d]^{2,q^4t^2}&\alpha Y_{1,1}Y_{2,q^5t^3}^{-1}\ar[dl]^{1,qt}   &Z_{1,t^4}^{-1} Z_{2,t^2} \ar[dd]^{2,t^4}
\\ \alpha Y_{1,q^2t^2}^{-1}Y_{2,qt}Y_{2,q^5t^3}^{-1}\ar[d]^{2,q^2t^2}&&
\\ Y_{2,q^3t^3}^{-1}Y_{2,q^5t^3}^{-1}&& Z_{2,t^6}^{-1}}$$
For the type $A_2$ with $r = 2$, we choose $\phi(1) = 0$, $\phi(2) = 1$ and we have the following.
$$\xymatrix{
Y_{1,1} \ar[d]^{1,q^2t}                        &Y_{1,1} \ar[d]^{1,q^2}                        &Z_{1,1} \ar[d]^{1,t}
\\Y_{1,q^4t^2}^{-1}Y_{2,q^2t}\ar[d]^{2,q^4t^2} &Y_{1,q^4}^{-1}Y_{2,q^2}\ar[d]^{2,q^4} &
Z_{1,t^2}^{-1}Z_{2,t}\ar[d]^{2,t^2}  
\\Y_{2,q^6t^3} & Y_{2,q^6}  & Z_{2,t^3}^{-1}
}$$
The following example was considered in \cite{fr2} (it is rewritten
here in the language of $q$-characters and twisted
$t$-characters). The type is $B_2^{(1)} = C_2^{(1)}$ and its Langlands
dual $D_3^{(2)} = A_3^{(2)}$. We have $\phi(1) = 0$, $\phi(2) =
1$. $\Pi_q$ gives the $q$-character of a fundamental $5$-dimensional
representation of $\U_q(C_2^{(1)})$ from \cite{fm} (see also
\cite{ks}) and $\Pi_t$ gives the following interpolating $(q,t)$-character.
$$\xymatrix{Y_{1,1} \ar[d]^{1,q^2t}                               &Y_{1,1} \ar[d]^{1,q^2}                             &
Z_{1,1} \ar[d]^{1,t}
\\Y_{1,q^4t^2}^{-1}Y_{2,qt}Y_{2,q^3t}\ar[d]^{2,q^4t^2}            &Y_{1,q^4}^{-1}Y_{2,q}Y_{2,q^3}\ar[d]^{2,q^4}       &
Z_{1,t^2}^{-1} Z_{2,t^2}\ar[dd]^{2,t^4}
\\\alpha Y_{2,qt}Y_{2,q^5t^3}^{-1}\ar[d]^{2,q^2t^2}           &Y_{2,q}Y_{2,q^5}^{-1}\ar[d]^{2,q^2}                &
\\Y_{2,q^3t^3}^{-1}Y_{2,q^5t^3}^{-1}Y_{1,q^2t^2}\ar[d]^{1,q^4t^3} &Y_{2,q^3}^{-1}Y_{2,q^5}^{-1}Y_{1,q^2}\ar[d]^{1,q^4}&
Z_{2,t^6}^{-1} Z_{1,-t^2}\ar[d]^{1,-t^3}
\\Y_{1,q^6t^4}^{-1}                                               &Y_{1,q^6}^{-1}                                     &
Z_{1,-t^4}^{-1}
}$$
By \cite{hgen} this is the twisted $t$-character of a fundamental
$4$-dimensional representation of $\U_t(A_3^{(2)})$.

Let us give another example for this type.
$$\xymatrix{Y_{2,1}Y_{2,q^2} \ar[d]^{2,q^3t} &&   
\\\alpha Y_{2,1}Y_{2,q^4t^2}^{-1}Y_{1,q^3t} \ar[d]^{2,qt}\ar[dr]^{1,q^5t^2}&& 
\\ Y_{2,q^2t^2}^{-1}Y_{2,q^4t^2}^{-1}Y_{1,qt}Y_{1,q^3t}\ar[d]^{1,q^3t^2}\ar[dr]^{1,q^5t^2}& \alpha  Y_{2,1}Y_{2,q^6t^2}Y_{1,q^7t^3}^{-1}\ar[d]^{2,qt}\ar[dr]^{2,q^7t^3} & 
\\ Y_{1,q^3t}Y_{1,q^5t^3}^{-1}\ar[d]^{1,q^5t^2} & Y_{2,q^2t^2}^{-1}Y_{2,q^6t^2}Y_{1,q^7t^3}^{-1}Y_{1,qt}\ar[dl]^{1,q^3t^2}\ar[d]^{2,q^7t^3} & \alpha Y_{2,1}Y_{2,q^8t^4}^{-1}\ar[dl]^{2,qt}
\\Y_{1,q^5t^3}^{-1}Y_{1,q^7t^3}^{-1}Y_{2,q^4t^2}Y_{2,q^6t^2}\ar[d]^{2,q^7t^3}& \alpha Y_{2,q^2t^2}^{-1}Y_{2,q^8t^4}^{-1}Y_{1,qt}\ar[dl]^{1,q^3t^2}&
\\\alpha Y_{1,q^5t^3}^{-1}Y_{2,q^4t^2}Y_{2,q^8t^4}^{-1}\ar[d]^{2,q^5t^3}&&
\\Y_{2,q^6t^4}^{-1}Y_{2,q^8t^4}^{-1}&&}$$
Here we have to check that it is in the $\mathfrak{K}$, since {\em a priori} it is
unclear that
$$\alpha  Y_{2,1}Y_{2,q^6t^2}Y_{1,q^7t^3}^{-1}
+ Y_{2,q^2t^2}^{-1}Y_{2,q^6t^2}Y_{1,q^7t^3}^{-1}Y_{1,qt}+ \alpha Y_{2,1}Y_{2,q^8t^4}^{-1}
+ \alpha Y_{2,q^2t^2}^{-1}Y_{2,q^8t^4}^{-1}Y_{1,qt}$$
is in $\mathfrak{K}_{2,q,t}$.
But if we subtract $\alpha Y_{2,1}(1 + A_{2,qt}^{-1})Y_{2,q^6t^2}(1 +
A_{2,q^7t^3}^{-1})Y_{1,q^7t^3}^{-1}\in\mathfrak{K}_{2,q,t}$, we get
$$(1 - \alpha) Y_{2,q^2t^2}^{-1}Y_{2,q^6t^2}Y_{1,q^7t^3}^{-1}Y_{1,qt}
= (1 - \alpha)Y_{1,q^7t^3}^{-1}Y_{1,qt}\in \mathfrak{K}_{2,q,t}.$$
$\Pi_q$ gives the $q$-character of a $11$-dimensional
KR module over $\U_q(B_2^{(1)})$ (it follows from
\cite{hcr} that the formula of \cite{kos, knh} is satisfied) and
$\Pi_t$ gives the following.
$$\xymatrix{Z_{2,-1} \ar[dd]^{2,-t^2} & 
\\ &
\\ Z_{2,-t^4}^{-1} Z_{1,-\epsilon t}Z_{1,\epsilon t}\ar[d]^{1,-\epsilon t^2}\ar[dr]^{1,\epsilon t^2}& 
\\ Z_{1,\epsilon t}Z_{1,-\epsilon t^3}^{-1}\ar[d]^{1,\epsilon t^2} &
Z_{1,\epsilon t^3}^{-1}Z_{1,-\epsilon t}\ar[dl]^{1,-\epsilon t^2}
\\ Z_{1,-\epsilon t^3}^{-1}Z_{1,\epsilon t^3}^{-1} Z_{2,-t^4}\ar[dd]^{2,-t^6}& 
\\&
\\Z_{2,-t^8}^{-1}&&}$$
By \cite{hgen} this is the twisted $t$-character of a fundamental $6$-dimensional
representation of $\U_t(A_3^{(2)})$.

\section{Interpolating $(q,t)$-characters}\label{interpolate}

In this section we construct interpolating $(q,t)$-characters in a
systematic way: we prove the existence and construct sums in $\mathfrak{K}$ 
with a unique dominant monomial which can be seen as
interpolating $(q,t)$-characters of virtual representations (Theorem
\ref{con}). Their existence implies Conjecture \ref{char} in many
cases. We will prove in Section \ref{endpf} that Theorem
\ref{con} implies Theorem \ref{thdual}.

Let us explain the main ideas of the construction of interpolating
$(q,t)$-characters. In \cite[Section 5]{hadv} a process is given to
construct some deformations of $q$-characters. Although the notion of
``interpolating $(q,t)$-characters'' considered in the present paper
is completely different from that of the ``$q,t$-characters'' in
\cite{hadv}, we use an analogous process (note that the
``$q,t$-characters'' of \cite{hadv} were first introduced in
\cite{Nab} for simply-laced affine quantum algebras by a different
method). In fact, the process of \cite{hadv} may be seen as a general
process to produce $t$-deformations under certain conditions. It is
based on an algorithm which is analogous to the Frenkel--Mukhin
algorithm for $q$-characters \cite{fm}.

Let us give the main points of the construction. We define a certain
property $P(n)$ depending on the rank $n$ of the Lie algebra which
means the existence of interpolating $(q,t)$-characters in $\mathfrak{K}$. 
To prove it by induction on $n$, assuming the
existence for the fundamental representations, we first construct some
elements $E(m)$ which are analogues of interpolating
$(q,t)$-characters for standard modules (tensor products of
fundamental representations). Then we have three additional steps:

{\it Step 1}: we prove $P(1)$ and $P(2)$ using a more precise property
$Q(n)$ such that $Q(n)$ implies $P(n)$. The property $Q(n)$ has the
following advantage: it can be checked by computation in elementary
cases $n=1,2$.

{\it Step 2}: we give some consequences of $P(n)$ which will be used
in the proof of $P(r)$ ($r>n$).

{\it Step 3}: we prove $P(n)$ ($n \geq 3$) assuming that $P(r)$,
$r\leq n$ are true. We give an algorithm to construct explicitly the
interpolating $(q,t)$-characters by using ideas of \cite{hadv}. As we
do not know {\it a priori} that the algorithm is well-defined in the
general case, we have to show that it never fails. This is a
consequence of $P(2)$ as it suffices to check the compatibility
conditions for pairs of nodes of the Dynkin diagram. Finally, we prove
that the algorithm stops, that is to say it gives a finite sum which
makes sense in $\mathfrak{K}$.

\subsection{Statement} In this section we prove, for $m\in B$, the
existence of an element $F(m)\in\mathfrak{K}$ such that $m$ is the
unique dominant monomial of $F(m)$. This will imply Theorem
\ref{thdual}. 

We have a partial ordering on the monomials of $\wt{\Yim}_{q,t}$:
$$m\leq m' \Leftrightarrow m(m')^{-1} \in \ZZ[\wt{A}_{i,a}^{-1},
  \alpha]_{i\in I, a\in\mathcal{C}}.$$

\begin{lem} A non-zero $\chi$ in $\mathfrak{K}_{i,q,t}$ has at
  least one $i$-dominant monomial.
\end{lem}

\demo Take a monomial $m$ in $\chi$ maximal for the partial ordering
$\leq$. It occurs in a product of generators of
$\mathfrak{K}_{i,q,t}$, whose product $M$ of highest monomials are
greater or equal to $m$ for the partial ordering, that is $Mm^{-1}$ is
a product of $v(M)$ factors $\wt{A}_{i,a}^{-1}$.  Let $N$ be the
maximal $v(M)$. We suppose that we have written $\chi$ so that $N$ is
minimal.  If $N = 0$, one of the products $M$ is equal to $m$, so $m$
is $i$-dominant. Otherwise, $N > 0$.  The products $M$ such that $v(M)
= N$ should cancel as $m$ is maximal in $\chi$. But the only case
where generators of $\mathfrak{K}_{i,q,t}$ have the same highest
monomial is when $i\in I_1$ as the dominant monomial $\alpha
Y_{i,a}Y_{i,aq^2}$ is the highest monomial of
$$\alpha Y_{i,a}(1 + \wt{A}_{i,aqt}^{-1}) \alpha Y_{i,aq^2}(1 +
\wt{A}_{i,aq^3t}^{-1})$$
and of
$$\alpha Y_{i,a}Y_{i,aq^2}(1 + \alpha \wt{A}_{i,aq^3t}^{-1} +
\wt{A}_{i,aq^3t}^{-1}\wt{A}_{i,aqt}^{-1}).$$ But the difference of the
two is $\alpha Y_{i,a}Y_{i,aq^2}\wt{A}_{i,aqt}^{-1} = \alpha
\prod_{j\sim i} Y_{j,aq}$ in $\tilde{\Yim}_{a,t}$. This monomial is
$i$-dominant in $\mathfrak{K}_{i,q,t}$ and strictly lower than $\alpha
Y_{i,a}Y_{i,aq^2}$. So we can rewrite the expression in such a way
that the new maximum of the $v(M)$ is strictly lower than
$N$. This is a contradiction.  \qed

For $J\subset I$, let $\Glie_J$ be the
semi-simple Lie algebra of Cartan Matrix $(C_{i,j})_{i,j\in J}$ and
$\U_q(\wh{\Glie})_J$ the associated quantum affine algebra with
coefficient $(r_i)_{i\in J}$.

As above, by considering a maximal monomial for the partial ordering,
we get the following:

\begin{lem}\label{leasto} A non-zero element of $\mathfrak{K}_J$ has
  at least one $J$-dominant monomial.
\end{lem}

For a monomial $m$ there is a finite number of monomial $m'\in
m\ZZ[\wt{A}_{i,a}^{-1}]_{a\in\mathcal{C}}$ which are $i$-dominant. 
Let $m = \prod_{i\in I,a\in\mathcal{C}}Y_{i,a}^{u_{i,a}(m)}$ and let
$C(m) = \{a\in\mathcal{C}|\exists i\in I,u_{i,a}(m)\neq 0\}$. Then we set
$$D(m) = \{m\wt{A}_{i_1,a_1}^{-1}\cdots \wt{A}_{i_N,a_N}^{-1}|N\geq 0,
i_j\in I, a_j\in C(m)q^{\NN^*}t^{\NN^*}\}.$$ 
Note that $D(m)$ is countable, any $m'\in D(m)$ satisfies $m'\leq m$
and $D(m') \subset D(m)$. Finally, set 
$$\wt{D}(m) = \oplus_{m'\in D(m)}\ZZ m'.$$ 
We prove the following result as in \cite[Lemma 3.14]{hadv}.

\begin{lem}\label{findom} For any monomial $m$, the set $D(m)\cap B$
  is finite.
\end{lem}

Let us state the main result of this section.

\begin{thm}\label{con} For all $n\geq 1$ we have the following
  property $P(n)$: for all semi-simple Lie algebras $\Glie$ of rank
  $\text{rk}(\Glie) = n$ and for all $m\in B$ there is a
  unique $F(m)\in\mathfrak{K}\cap \tilde{D}(m)$ such that $m$ is the
  unique dominant monomial of $F(m)$.
\end{thm}

\begin{rem} If $m$ is of the form $\alpha m'$, then the existence of
$F(m)$ follows from the analogous result for the
$q$-characters. Indeed, in \cite{hadv} an algorithm inspired by the
Frenkel--Mukhin algorithm \cite{fm} was proposed (as well as its
$t$-deformation in the sense of \cite{hadv}): if it is well-defined,
then for a dominant monomial $m\in\ZZ[Y_{i,q^r}]_{i\in I,r\in\ZZ}$ it
gives $F(m)$ in the ring of $q$-character such that $m$ is the unique
dominant monomial of $F(m)$ (see also \cite{haja}). As a consequence,
it suffices to prove the result when $m$ is a product of the
$W_{i,a}$.
\end{rem}

\subsection{Proof of Theorem \ref{con}}

First note that for $n=1$ we have already proved this result. For a
general $n$, the uniqueness follows from lemma \ref{leasto}.

First, we define a new property $Q(n)$.

\begin{defi} For $n\geq 1$ denote by $Q(n)$ the property ``for all
  semi-simple Lie algebras $\Glie$ of rank $n$, for all $i\in I$ there
  is a unique $F(W_{i,1})\in\mathfrak{K}\cap
  \wt{D}(W_{i,1})$ such that $W_{i,1}$ is the unique dominant monomial
  of $F(W_{i,1})$.''\end{defi}

\subsubsection{Construction of the $E(m)$}\label{conste}

We suppose that for $i\in I$, there is $F(W_{i,1})\in\mathfrak{K} \cap
\wt{D}(W_{i,1})$ such that $W_{i,1}$ is the unique dominant monomial
of $F(W_{i,1})$ (that is the property $Q(n)$ is satisfied).

For $a\in\mathcal{C}$ consider $s_a:\wt{\Yim}_{q,t}\rightarrow
\wt{\Yim}_{q,t}$ the algebra morphism such that $s_a(Y_{j,b}) =
Y_{j,ab}$. We can define for $m = \prod_{i\in
I,a\in\mathcal{C}}W_{i,a}^{w_{i,a}}$ the element
$$E(m) = \prod_{i\in I, a\in\mathcal{C}} (s_a(F(W_{i,1})))^{w_{i,a}}
\in \mathfrak{K}\cap(\prod_{i\in I,a\in\mathcal{C}}(\wt{D}(W_{i,a}))^{w_{i,a}})
\subset \mathfrak{K}\cap \wt{D}(m).$$

\subsubsection{Step 1} 

First, we prove that $Q(n)$ implies $P(n)$.

\begin{lem} For $n\geq 1$, property $Q(n)$ implies property
  $P(n)$.\end{lem}

\demo We suppose that $Q(n)$ is true. In particular, we can construct
$E(m)\in\mathfrak{K} \cap \wt{D}(m)$ for $m\in B$ as above. Let us
prove $P(n)$. Let $m\in B$. The uniqueness of $F(m)$ follows from
Lemma \ref{leasto}. Let $m_L=m>m_{L-1}>\cdots >m_1$ be the dominant
monomials of $D(m)$ with a total ordering compatible with the partial
ordering (it follows from Lemma \ref{findom} that $D(m)\cap B$ is
finite). Let us prove by induction on $l$ the existence of
$F(m_l)$. The unique dominant monomial of $D(m_1)$ is $m_1$, so
$F(m_1)=E(m_1)\in \wt{D}(m_1)$. In general, let $\lambda_1,\cdots
,\lambda_{l-1}\in\ZZ$ be the coefficient of the dominant monomials
$m_1,\cdots ,m_{l-1}$ in $E(m_l)$. We put
$$F(m_l)=E(m_l)-\underset{r=1\cdots l-1}{\sum}\lambda_r F(m_r).$$ It
follows from the construction that $F(m)\in \wt{D}(m)$ because for
$m'\in D(m)$ we have $E(m')\in \wt{D}(m')\subseteq \wt{D}(m)$.  \qed

\begin{cor} The properties $Q(1)$, $Q(2)$, and hence $P(1)$, $P(2)$, are
  true.\end{cor}

\noindent This allow us to start our induction in the proof of Theorem
\ref{con}.

\demo For $n = 1$ we have two cases $A_1$ with $r = 1$ or $r = 2$. The
explicit formulas have been given above. For $n = 2$ we have five cases
$A_1\times A_1$ with $r = 1,2$, $A_2$ with $r = 1, 2)$, $B_2$. The
cases $A_1\times A_1$ are a direct consequence of the case $n =
1$. For $A_2$, $i = 1, 2$ are symmetric so it suffices to give the
formulas for $i = 1$ as we did above. We also gave the formulas for
$B_2$ above.\qed

\subsubsection{Step 2}\label{consp} Let be $n\geq 1$. We suppose in
this section that $P(n)$ is proved. We give some consequences of
$P(n)$ which will be used in the proof of $P(r)$ ($r>n$).

From Lemma \ref{findom}, an element of $\wt{\Yim}_{q,t}$ has
a finite number of dominant monomials.

\begin{prop}\label{thth} We suppose $\text{rk}(\Glie) = n$. We have
$$\mathfrak{K} = {\bigoplus}_{m\in B}\ZZ F(m).$$
\end{prop}

\demo Let $\chi\in \mathfrak{K}$. Let $m_1,\cdots ,m_L\in B$ the
dominant monomials occurring in $\chi$ and
$\lambda_1,\cdots,\lambda_L\in\ZZ$ their coefficients. It follows from
Lemma \ref{leasto} that $\chi = \underset{l=1\cdots L}{\sum}\lambda_l
F(m_l)$.  \qed

\begin{cor}\label{recufond} We suppose $|I|> n $ and let $J\subset I$
  such that $|J|=n$. For $m\in B_J$, there is a unique
  $F_J(m)\in\mathfrak{K}_J$ such that $m$ is the unique $J$-dominant
  monomial of $F_J(m)$. Moreover, $F_J(m)\in \wt{D}(m)$ and we have 
$$\mathfrak{K}_J = \bigoplus_{m\in B_J}\ZZ F_J(m).$$
\end{cor}

\demo The uniqueness of $F(m)$ follows from lemma \ref{leasto}. Let us
write $m = m_Jm'$ where $m_J=\underset{i\in
J,l\in\ZZ}{\prod}Y_{i,l}^{u_{i,l}(m)}\in B_J$. In
particular, Proposition \ref{thth} with the algebra
$\U_q(\wh{\Glie})_J$ of rank $n$ gives $m_J\chi$, where $\chi$ is a
polynomial in the variable $\wt{A}_{i,l}^{-1}$ for
$\U_q(\wh{\Glie})_J$. It suffices to put $F_J(m)=m\nu_J(\chi)$, where
$\nu_J$ is the ring morphism which sends a variable $\wt{A}_{i,a}^{-1}$
for $\U_q(\wh{\Glie})_J$ to the corresponding variable for
$\U_q(\wh{\Glie})$. The last assertion is proved as in Proposition
\ref{thth}.\qed

\subsubsection{Step 3} We explain why properties $P(r)$ ($r<n$) imply
$P(n)$. In particular, we define an algorithm which constructs
explicitly the $F(m)$ by using ideas of \cite{hadv}.

We prove the property $P(n)$ by induction on $n\geq 1$. We have proved
$P(1)$ and $P(2)$. Let $n\geq 3$ and suppose that $P(r)$ is proved for
$r < n$.

Let $m_0 \in B$ and $m_0, m_1, m_2,\cdots $ the countable set
$D(m_0)$ with indexes such that $m_j\geq m_{j'}$ implies $j'\geq j$.

For $J\varsubsetneq I$ and $m\in B_J$, it follows from $P(r)$ and
corollary \ref{recufond} that there is a unique $F_J(m)\in
\wt{D}(m)\cap\mathfrak{K}_J$ such that $m$ is the unique $J$-dominant
monomial of $F_J(m)$ and that $\mathfrak{K}_J = \bigoplus_{m\in
B_J}\ZZ F_J(m)$.  If $m\notin B_J$, we denote $F_J(m)=0$. For
$\chi\in\wt{\Yim}_{q,t}$, $[\chi]_{m'}\in\ZZ$ is the coefficient of
$m'$ in $\chi$.

We consider the following inductive definition of the sequences
$(s(m_r))_{r\geq 0}\in\ZZ^{\NN}$, $(s_J(m_r))_{r\geq 0}\in\ZZ^{\NN}$
($J\varsubsetneq I$),
$$s(m_0) = 1\text{ , }s_J(m_0) = 0,$$
and for $r\geq 1$, $J\subsetneq I$,
$$s_J(m_r) =
\underset{r'<r}{\sum}(s(m_{r'})-s_J(m_{r'}))[F_J(m_{r'})]_{m_r},$$
$$s(m_r) = \begin{cases} s_J(m_r)&\text{ if $m_r\notin B_J$,}
\\ 0&\text{ if $m_r\in B$.}
\end{cases}$$
The definition of $s_J$ means that we add the various contributions of
the $m_{r'}$ where $r' < r$ with coefficient $(s(m_{r'})-s_J(m_{r'}))$,
so that a contribution is not counted twice. For the definition of $s(m_r)$,
there is something to be proved, that it that the various $s_J(m_r)$ for
$m_r\notin B_J$ coincide.

We prove that the algorithm defines sequences in a unique way. We see
that if $s(m_r)$, $s_J(m_r)$ are defined for $r\leq R$, then so are
$s_J(m_{R+1})$ for $J\subsetneq I$. Moreover, $s_J(m_R)$ imposes the
value of $s(m_{R+1})$, and by induction the uniqueness is clear. We
say that the algorithm is well-defined to step $R$ if there exist
$s(m_{r})$, $s_J(m_r)$ such that the formulas of the algorithm are
satisfied for $r\leq R$.

\begin{lem} The algorithm is well-defined to step $r$ if and only if
$$\forall J_1,J_2\varsubsetneq I, \forall r'\leq r, (m_{r'}\notin
B_{J_1} \text{ and }m_{r'}\notin B_{J_2}\Rightarrow s_{J_1}(m_{r'}) =
s_{J_2}(m_{r'})).$$\end{lem}

\demo If for $r'<r$ the $s(m_{r'})$, $s_J(m_{r'})$ are well-defined,
so is $s_J(m_r)$. If $m_r\in B$, $s(m_r) = 0$ is well-defined. If
$m_r\notin B$, it is well-defined if and only if $\{s_J(m_r)
|m_r\notin B_J\}$ has a unique element.\qed

If the algorithm is well-defined to step $r$, then for $J\varsubsetneq
I$ we set
$$\mu_J(m_r) = s(m_r) - s_J(m_r)\text{ , }\chi_J^r=\underset{r'\leq
r}{\sum}\mu_{J}(m_{r'})F_J(m_{r'})\in\mathfrak{K}_J.$$

We prove as in \cite[Lemma 5.21]{hadv} (except that the coefficients
are in $\ZZ$ and not in $\ZZ[t^{\pm 1}]$) the following:

\begin{lem} If the algorithm is well-defined to step $r$, for
  $J\subset I$ we have
$$\chi_J^r\in (\underset{r'\leq r}{\sum}
s(m_{r'})m_{r'})+s_J(m_{r+1})m_{r+1}+\underset{r'>r+1}{\sum}\ZZ
m_{r'}$$ For $J_1\subset J_2\subsetneq I$, we have
$$\chi_{J_2}^r=\chi_{J_1}^r+\underset{r'>r}{\sum}\lambda_{r'}
F_{J_1}(m_{r'})$$ where $\lambda_{r'}\in\ZZ$. In particular, if
$m_{r+1}\notin B_{J_1}$, we have $s_{J_1}(m_{r+1})=s_{J_2}(m_{r+1})$.
\end{lem}

We prove as in \cite[Lemma 5.22]{hadv} the following.

\begin{lem}\label{nfail} The algorithm never fails.\end{lem}

Now we aim at proving that the algorithm stops. We will use the
following notion \cite{fm}:

\begin{defi}\label{monomrn} A non-trivial $m = \prod_{i\in
    I,a\in\CC^\times}Y_{i,a}^{u_{i,a}(m)}$ is said to be right-negative if
  for all $a\in\CC^\times, j\in I$ we have $(u_{j,aq^{L_a}}(m)\neq
  0\Rightarrow u_{j,aq^{L_a}}(m)<0)$ where
$$L_a=\text{max}\{l\in\ZZ|\exists i\in I, u_{i,aq^L}(m)\neq
0\}.$$\end{defi} 

$D(m)$ is graded by finite-dimensional subspaces such that the degree
of the monomial $m' = m\wt{A}_{i_1,a_1}^{-1}\cdots
\wt{A}_{i_N,a_N}^{-1}$ in $D(m)$ is $N$. Then we can consider the
corresponding graded completion $\overline{D}(m)$ of $\wt{D}(m)$. By
an infinite sum in $\tilde{\Yim}_{q,t}$ we mean an element in such a
completion. We have analogous definitions for infinite sums in
$\Yim_q$ and in $\Yim_t^L$

\begin{lem}\label{finites} Let $S$ be an infinite sum in $\Yim_q$
  (resp. in $\Yim_t^L$) which is an infinite sum of
  elements in $\mathfrak{K}_{i,q}$ (resp. in
  $\mathfrak{K}_{i,t}^L$) for any $i\in I$. If $S$ contains a finite number of dominant
  monomials, then $S$ is a finite sum in $\Yim_q$ (resp. in
  $\Yim_t^L$).
\end{lem}

\demo We prove the result for $\Yim_q$ (the proof is completely analogous
for $\Yim_t^L$ by using results in \cite{hgen}). Let $m_1,\cdots,m_L$
be the dominant monomials occurring in $S$ and
$\lambda_1,\cdots,\lambda_L$ their multiplicity. For $m$ a dominant
monomial, there is $F_q(m)\in\text{Im}(\chi_q)$ with a unique dominant
monomial $m$ (see the construction in \cite[Section 5.1]{hadv} by
using $q$-characters which are finite sums). Then 
$$S' = S -
\sum_{1\leq l\leq L}\lambda_l F_q(m_l)$$ 
has no dominant monomial and
for any $i\in I$ is an infinite sum of elements in
$\mathfrak{K}_{i,q}$. So if $S'\neq 0$, a maximal monomial occurring
in $S'$ is dominant, contradiction. So $S' = 0$.\qed

Now we can prove the following:

\begin{lem} The algorithm stops and $\chi=\underset{r\geq
    0}{\sum}s(m_r)m_r\in\mathfrak{K}\cap \wt{D}(m_0)$. Moreover, the
  only dominant monomial in $\chi$ is $m_0$.
\end{lem}

\demo Consider the ({\em a priori}, non necessarily finite) sum $\chi$ in
$\overline{D}(m_0)$. We prove as in \cite[Lemma 5.23]{hadv} that for
each $i\in I$, $\chi$ is an infinite sum of elements in
$\mathfrak{K}_{i,q,t}$.

There in $N\in\ZZ$ such that $m_0\in \ZZ[Y_{i,q^r t^l}]_{i\in I,
r,l\leq N }$.  By construction with the algorithm, only a finite
number of monomials of $F(m_0)$ are in
$m_0\ZZ[\tilde{A}_{i,q^rt^l}^{-1}]_{r\leq N\text{ or }l\leq N}$.  Let
us consider another monomials $m'\notin
\ZZ[\tilde{A}_{i,q^rt^l}^{-1}]_{r\leq N\text{ or }l\leq N}$ occurring
in $\chi$. The specializations $\Pi_q(m')$ and $\Pi_t(m')$ are
right-negative. Indeed for any $r_1, r_2 > N$ and $j\in I$, the
specializations of $m_0\wt{A}_{j,q^{r_1}t^{r_2}}^{-1}$ are
right-negative. Moreover the specializations of the
$\wt{A}_{i,a}^{-1}$ are right-negative, and a product of
right-negative monomials is right-negative \cite{fm}. Since a
right-negative monomial is not dominant, we can conclude that the
specializations of $m'$ are not dominant.  So $\Pi_q(\chi)$ and
$\Pi_t(\chi)$ have a finite number of dominant monomials.  So these
are finite sums by Lemma \ref{finites}.  As $\tilde{\Yim}_{q,t}$ is
obtained by a quotient by $\text{Ker}(\Pi_q)\cap \text{Ker}(\Pi_t)$,
$\chi$ is a finite sum.  \qed

This lemma implies the following.

\begin{cor} For $n\geq 3$, if the $P(r)$ ($r<n$) are true, then $P(n)$
  is true.\end{cor}
In particular, Theorem \ref{con} is proved by induction on $n$. 

\subsection{Proof of Theorem \ref{thdual}}\label{endpf}

Let us explain how Theorem \ref{con} implies Theorem
\ref{thdual}. 

First consider the dominant interpolating monomial
$$m = W_{i,a}W_{i,at^2q^4}\cdots W_{i,a(t^2q^4)^{k-1}}.$$ 
The specializations by $\Pi_q,\Pi_t$ of $m$ correspond to the highest
monomials of KR modules respectively over
$\U_q(\wh{\Glie})$ and $\U_t({}^L\wh{\Glie})$. 
By construction, the monomials $m'$ occurring in $F(m) - m$ are of the form
$$m' = (m\wt{A}_{i,aq_i^{2k}t}^{-1})\wt{A}_{i_1,a_1}^{-1}\cdots
\wt{A}_{i_N,a_N}^{-1} \text{ where $i_1,\cdots,i_N\in I$ and
$a_1,\cdots,a_N\in\mathcal{C}$.}$$ 
As a consequence, $\Pi_q(m')$ and
$\Pi_t(m')$ are right-negative.  Indeed, the specialization of
$m\wt{A}_{i,aq_i^{2k - 1}t}^{-1}$ and of the $\wt{A}_{i,a}^{-1}$ are
right-negative, and a product of right-negative monomials is
right-negative \cite{fm}. Since a right-negative monomial is not
dominant, the specializations of $F(m)$ are affine-minuscule. 
By Theorem \ref{krun} and Theorem \ref{krdeux}, this
completes the proof of the first statement of Theorem \ref{thdual}
for KR modules.

Now we have the following compatibility property with
tensor products.

\begin{prop}\label{tens} Let $V_1$ and $V_2$ be two simple representations of
$\U_q(\hat{\Glie})$ with respective Langlands dual representations
$V_1^L$ and $V_2^L$. 
If $V_1\otimes V_2$ is simple, then $V_1^L\otimes V_2^M$ is a
Langlands dual representation to $V_1\otimes V_2$.
\end{prop}

\demo For $\chi_1$ and $\chi_2$ interpolating $(q,t)$-characters
respectively of $V_1$ and $V_2$, the product $\chi_1\chi_2$ is
clearly an interpolating $(q,t)$-character of $V_1\otimes V_2$ as
$\mathfrak{K}$ is a subring of $\tilde{\Yim}_{q,t}$ and $\Pi_q$ is a
ring morphism. We can conclude for the last point as $\Pi_t$ is a ring
morphism.  \qed

According to the above discussion, 
this completes the proof of the first statement of Theorem \ref{thdual}.

\begin{rem}\label{uni} If $V_1^L$ and $V_2^L$
are non trivial, then the uniqueness statement does not hold
for $V_1\otimes V_2$: by shifting the spectral parameter in $\chi_1$
by $t^N$ where $N\neq 1$ without changing $\chi_2$, we get another
Langlands dual representation which can not obtained from $V_1^L\otimes V_2^L$
by a shift of spectral parameter.
\end{rem}

Now consider a KR module $V$ of highest monomial
$W_{i,a}W_{i,aq^4}\cdots W_{i,a(q^4)^{k-1}}$.
Suppose that we have a Langlands dual representation
$V^L$ which is not a KR module and consider a
corresponding interpolating $(q,t)$-character $\chi$.
Let $M$ be the highest monomial of $\chi$. 
Then there are
$b\neq b'$ such that $W_{i,b}$, $W_{i,b'}$ occur in $M$
but $W_{i,bq^4t^2}$, $W_{i,b'q^4t^2}$ do not occur in $M$.
Suppose that $r_i = 1$ (resp. $r_i = 2$).
As a consequence, $MA_{i,bq^3t}^{-1}A_{i,bqt}^{-1}$
and $MA_{i,b'q^3t}^{-1}A_{i,b'qt}^{-1}$ 
(resp. $MA_{i,bq^2t}^{-1}$ and $MA_{i,b'q^2t}^{-1}$)
occur in $\chi$ and have distinct image by $\Pi_q$.
But by \cite[Lemma 5.5]{hcr}, $\chi_q(V)$ contains a unique monomial
of the form $\Pi_q(M)A_{i,a'}^{-1}A_{i,a''}^{-1}$ (resp. $\Pi_q(M)A_{i,a'}^{-1}$).
This contradicts $\Pi_q(\chi) = \chi_q(V)$.

\subsection{Additional comments}\label{comments}
Note that it is easy to construct interpolating $(q,t)$-characters of
non-simple representations by using tensor products of
KR modules which are not simple, by the same method as in 
Proposition \ref{tens}. More interestingly,
to illustrate Conjecture \ref{char}, let us give an example of a
simple non affine-minuscule module which satisfies the Langlands
duality. Consider the $\U_q(C_2^{(1)})$-module $V =
L(Y_{1,1}^2Y_{2,q^5}Y_{2,q^7})$. Note that
$L(Y_{2,1}Y_{2,q^2}Y_{1,q^7})\otimes L(Y_{1,q^7})$ is simple as it is
affine-minuscule. So by \cite[Lemma 4.10]{hcmp}, $V \simeq
L(Y_{1,1}Y_{2,q^5}Y_{2,q^7})\otimes L(Y_{1,1})$. Moreover
$$\chi_q(L(Y_{1,1})\otimes L(Y_{2,q^5}Y_{2,q^7})) =
\chi_q(L(Y_{2,1})\otimes L(Y_{2,7})) +
\chi_q(L(Y_{1,1}Y_{2,q^5}Y_{2,q^7})),$$
has they have the same multiplicity $1$ on the dominant monomials. So 
$$\text{dim}(L(Y_{2,q^5}Y_{2,q^7})) = 55 - 16 = 39\text{ and }\text{dim}(V) = 39 \times 5 = 195.$$
Now consider the $\U_t(A_3^{(2)})$ simple module $V^L = L(Z_{1,1}^2Z_{2,t^6})$. In the same way, by \cite[Proposition 4.7]{hprod}, we have $V^L\simeq L(Z_{1,1}Z_{2,t^6})\otimes L(Z_{1,1})$, and as
$$\chi_t^\sigma(L(Z_{1,1})\otimes L(Z_{2,t^6})) = \chi_t^\sigma(L(Z_{1,1}Z_{2,t^6})) - \chi_t^\sigma(L(Z_{1,-t^2})),$$
we get $\text{dim}(L(Z_{1,1}Z_{2,t^6})) = 24 - 4 = 20$ and $\text{dim}(V^L) = 20\times 4 = 80$. As for their dimension above, it is easy to compute the $q$-character (resp. twisted $t$-character) of $V$ (resp. $V^L$), and so to check that $V$ satisfies the Langlands duality with the Langlands dual module $V^L$. We do not list the $195$ monomial of the interpolating $(q,t)$-character, but the $80$ monomials which do not have $\alpha$ in their coefficient. It suffices to multiply one of the $4$ monomials of the sum
$$Y_{1,1} + Y_{1,q^4t^2}^{-1}Y_{2,qt}Y_{2,q^3t} + Y_{2,q^3t^3}^{-1}Y_{2,q^5t^3}^{-1}Y_{1,q^2t^2} + Y_{1,q^6t^4}^{-1}$$ by one the following $20$ monomials. We use the notation $i_a = Y_{i,a}$ (analog notation will also be used in the following).

$ 1_1 2_{q^5t^3} 2_{q^7t^3} , 
1_1 1_{q^6t^4} 1_{q^8t^4} 2_{q^7t^5}^{-1} 2_{q^9t^5}^{-1} ,
1_{q^4t^2}^{-1} 2_{qt} 2_{q^3t} 2_{q^5t^3} 2_{q^7t^3} 
1_1 1_{q^6t^4} 1_{q^{12}t^6}^{-1} 2_{q^7t^5}^{-1} 2_{q^{11}t^5} ,
\\1_1 1_{q^{10}t^6}^{-1} 1_{q^8t^4} ,
1_{q^4t^2}^{-1} 1_{q^6t^4} 1_{q^8t^4} 2_{qt} 2_{q^3t} 2_{q^7t^5}^{-1} 2_{q^9t^5}^{-1} ,
1_1 1_{q^{10}t^6}^{-1} 1_{q^{12}t^6}^{-1} 2_{q^9t^5} 2_{q^{11}t^5} ,
1_{q^4t^2}^{-1} 1_{q^{10}t^6}^{-1} 2_{qt} 2_{q^3t}1_{q^8t^4},
\\1_{q^4t^2}^{-1} 1_{q^6t^4} 1_{q^{12}t^6}^{-1} 2_{qt}2_{q^3t}2_{q^7t^5}^{-1}2_{q^{11}t^5},
1_{q^4t^2}^{-1}1_{q^{10}t^6}^{-1}1_{q^{12}t^6}^{-1}2_{qt}2_{q^3t}2_{q^9t^5}2_{q^{11}t^5},
\\1_{q^2t^2}1_{q^6t^4}1_{q^8t^4}2_{q^3t^3}^{-1}2_{q^5t^3}^{-1}2_{q^7t^5}^{-1}2_{q^9t^5}^{-1},
1_12_{q^{11}t^7}^{-1}2_{q^{13}t^7}^{-1},
1_{q^2t^2}1_{q^{10}t^6}^{-1}2_{q^3t^3}^{-1}2_{q^5t^3}^{-1}1_{q^8t^4},
\\1_{q^2t^2}1_{q^6t^4}1_{q^{12}t^6}^{-1}2_{q^3t^3}^{-1}2_{q^5t^3}^{-1}2_{q^7t^5}^{-1}2_{q^{11}t^5},
1_{q^4t^2}^{-1}2_{qt}2_{q^3t}2_{q^{11}t^7}^{-1}2_{q^{13}t^7}^{-1},
1_{q^2t^2}1_{q^{10}t^6}^{-1}1_{q^{12}t^6}^{-1}2_{q^3t^3}^{-1}2_{q^5t^3}^{-1}2_{q^9t^5}2_{q^{11}t^5},
\\1_{q^2t^2}2_{q^3t^3}^{-1}2_{q^5t^3}^{-1}2_{q^{11}t^7}^{-1}2_{q^{13}t^7}^{-1},
1_{q^6t^4}^{-1}1_{q^{10}t^6}^{-1}1_{q^8t^4},
1_{q^6t^4}^{-1}1_{q^{10}t^6}^{-1}1_{q^{12}t^6}^{-1}2_{q^9t^5}2_{q^{11}t^6},
1_{q^6t^4}^{-1}2_{q^{11}t^7}^{-1}2_{q^{13}t^7}^{-1}.
$

Nakajima \cite{Nab} has computed the $q$-characters of simple modules
from those of standard modules (tensor products of fundamental
representations) using quiver varieties. His results are not available
for non-simply laced untwisted quantum affine algebras, but the second
author has conjectured \cite{hadv} that analogous result do hold in
this case. The result on tensor products in Proposition \ref{tens} 
is an indication of the compatibility of the two conjectures.

\subsection{Example}\label{hgd} We give an example of an interpolating
$(q,t)$-character that we get for a Lie algebra or rank
strictly greater than $2$ by the process described in the
proof. Consider $\U_q(C_3^{(1)})$ with $\phi(1) = \phi(3) = 0$, $\phi(2) = 1$.
$$\xymatrix{Y_{3,1} \ar[d]^{3,q^2t}&
\\ Y_{3,q^4t^2}^{-1}Y_{2,qt}Y_{2,q^3t}\ar[d]^{2,q^4t^2}&
\\ \alpha Y_{1,q^4t^2}Y_{2,qt}Y_{2,q^5t^2}^{-1}\ar[d]^{2,q^2t^2}\ar[dr]^{1,q^5t^3}&
\\Y_{3,q^2t^2}Y_{2,q^3t^3}^{-1}Y_{2,q^5t^2}^{-1}Y_{1,q^4t^2}Y_{1,q^2t^2}\ar[d]^{3,q^4t^3}\ar[dr]^{1,q^5t^3}&\alpha^2Y_{1,q^6t^4}^{-1}Y_{2,qt}\ar[d]^{2,q^2t^2}
\\Y_{3,q^6t^4}^{-1}Y_{1,q^4t^2}Y_{1,q^2t^2}\ar[d]^{1,q^5t^3}&\alpha Y_{3,q^2t^2}Y_{2,q^3t^3}^{-1}Y_{1,q^6t^4}^{-1}Y_{1,q^2t^2}\ar[dl]^{3,q^4t^3}\ar[d]^{1,q^3t^3}
\\\alpha Y_{3,q^6t^4}^{-1}Y_{1,q^6t^4}^{-1}Y_{1,q^2t^2}Y_{2,q^5t^3}\ar[d]^{2,q^6t^4}\ar[dr]^{1,q^3t^3}& Y_{3,q^2t^2}Y_{1,q^6t^4}^{-1}Y_{1,q^4t^4}^{-1}\ar[d]^{3,q^4t^3}
\\\alpha^2 Y_{1,q^2t^2}Y_{2,q^7t^4}^{-1}\ar[dr]^{1,q^3t^3}& Y_{3,q^6t^4}^{-1}Y_{2,q^3t^3}Y_{2,q^5t^3}Y_{1,q^6t^4}^{-1}Y_{1,q^4t^4}^{-1}\ar[d]^{2,q^6t^4}
\\&\alpha Y_{2,q^3t^3}Y_{2,q^7t^5}^{-1}Y_{1,q^4t^4}^{-1}\ar[d]^{2,q^4t^4}
\\&Y_{2,q^5t^5}^{-1}Y_{2,q^7t^5}^{-1}Y_{3,q^4t^4}\ar[d]^{3,q^6t^5}
\\&Y_{3,q^8t^6}^{-1}}$$
The specialization at $t=1$ gives the $q$-character of a $14$-dimensional fundamental representation of $\U_q(C_3^{(1)})$ from
\cite{fm} (see also \cite{ks}). The specialization at $q = \epsilon$
gives the twisted $t$-character of a $8$-dimensional fundamental
representation of $\U_t(D_5^{(2)})$ \cite{hgen}.
$$\xymatrix{Z_{3,1} \ar[d]^{3,t}&
\\ Z_{3,t^2}^{-1}Z_{2,t^2}\ar[dd]^{2,t^4}&
\\ &
\\Z_{3,-t^2}Z_{2,t^6}^{-1}Z_{1,t^4} \ar[d]^{3,-t^3}\ar[ddr]^{1,t^6}&
\\Z_{3,-t^4}^{-1}Z_{1,t^4}\ar[ddr]^{1,t^6}&
\\& Z_{3,-t^2}Z_{1,t^8}^{-1}\ar[d]^{3,-t^3}
\\& Z_{3,-t^4}^{-1}Z_{2,t^6}Z_{1,t^8}^{-1}\ar[dd]^{2,t^8}
\\&
\\&Z_{2,t^{10}}^{-1}Z_{3,t^4}\ar[d]^{3,t^5}
\\&Z_{3,t^6}^{-1}}.$$

\section{Triple-laced case}\label{typeg}

Now we suppose that $r = 3$, that is to say we consider
$\U_q(G_2^{(1)})$ and its Langlands dual $\U_t(D_4^{(3)})$.  The
results and their proofs are completely analogous to the case $r = 2$,
except that we have to change some definitions and formulas and we
have to check the existence of interpolating $(q,t)$-characters in
some examples as we did for $r = 2$.

\subsection{Definitions of interpolating structures}

We set $\epsilon = e^{i\pi/3}$. For the Dynkin diagram of $G_2$ we use
the convention $r_1 = 3$ and $r_2 = 1$. We have $r_1^\vee = 1$ and
$r_2^\vee = 3$.

For the $q$-characters of $\U_q(G_2^{(1)})$ we have
$$A_{1,a} = Y_{1,aq^{-3}}Y_{1,aq^3}
Y_{2,aq^{-2}}^{-1}Y_{2,a}^{-1}Y_{2,aq^2}^{-1}\text{ , }A_{2,a} =
Y_{2,aq^{-1}}Y_{2,aq}Y_{1,a}^{-1},$$
$$\mathfrak{K}_{1,q} = \ZZ[Y_{1,a}(1 + A_{1,aq^3}^{-1}), Y_{2,a}^{\pm
    1}]_{a\in q^\ZZ}\text{ , } \mathfrak{K}_{2,q} = \ZZ[Y_{2,a}(1 +
    A_{2,aq}^{-1}), Y_{1,a}^{\pm 1}]_{a\in q^\ZZ}.$$ For the
    twisted $t$-characters of $\U_t(D_4^{(3)})$ we have
$$B_{1,a} = Z_{1,at^{-1}}Z_{1,at} Z_{2,a^3}^{-1}\text{ , }B_{2,a^3} =
Z_{2,a^3t^3}Z_{2,a^3t^{-3}}Z_{1,a}^{-1}
Z_{1,a\epsilon^2}^{-1}Z_{1,a\epsilon^4}^{-1},$$
$$\mathfrak{K}_{1,t}^L = \ZZ[Z_{1,a}(1 + Z_{1,at}^{-1}), Z_{2,a}^{\pm
1}]_{a\in\epsilon^\ZZ t^\ZZ}\text{ , }\mathfrak{K}_{2,t}^L =
\ZZ[Z_{2,a}(1 + A_{2,at^3}^{-1}), Z_{1,a}^{\pm 1}]_{a\in\epsilon^\ZZ
t^\ZZ}.$$ For $a\in\mathcal{C}$ let $W_{1,a} = Y_{1,a}$, $W_{2,a} =
Y_{2,aq^{-2}}Y_{2,a}Y_{2,a q^2}$.

Let us consider an interpolating map $\beta(q,t)$ such that
$\beta(q,1) = 1$ and $\beta(\epsilon, t) = 0$. We can use, for
example, the following map introduced in \cite{bpil}:
$$\beta(q,t) = \frac{(q^3 - q^{-3})(qt^{-1} -
q^{-1}t)(q^5t^{-1}-q^{-5}t)(q^4t^{-2} - q^{-4}t^2)}{(q -
q^{-1})(q^3t^{-1} - q^{-3}t)(q^4t^{-1} - tq^{-4})(q^5t^{-2} -
q^{-5}t^2)}.$$ Consider
$$\Yim_{q,t} = \ZZ[W_{i,a}^{\pm 1}, \beta Y_{i,a}^{\pm 1},
\beta]_{i\in I,a\in\mathcal{C}}.$$ 
We have the specializations maps $\Pi_q, \Pi_t$ and the ideal
$\text{Ker}(\Pi_q)\cap
\text{Ker}(\Pi_t)$ is generated by the elements
$$\beta(\beta - 1)\text{ , }\beta (Y_{i,a} -
Y_{i,at})\text{ , } (\beta - 1)(W_{i,aq} - W_{i,a\epsilon})\text{ , }
(W_{i,a} - W_{i,at})(W_{j,bq} - W_{j,b\epsilon}),$$
for $i,j\in I$ and $a,b\in\mathcal{C}$. We work in the ring
$\wt{\Yim}_{q,t} = \Yim_{q,t}/\langle\text{Ker}(\Pi_q)\cap
\text{Ker}(\Pi_t)\rangle$.

\begin{defi} We define for $a\in\mathcal{C}$ the interpolating root
  monomials
$$\wt{A}_{1,a} = Y_{i,a(q^3t)^{-1}}Y_{i,aq^3t} 
(Y_{2,aq^{-2}}Y_{2,a}Y_{2,aq^2})^{-1}\text{ , } \wt{A}_{2,a} =
Y_{2,a(qt)^{-1}}Y_{2,aqt}  Y_{1,a}^{-1}.$$
\end{defi}

We will use the identification $Z_{1,a} = Y_{1,a}$ and
$Y_{2,a}Y_{2,\epsilon^2 a}Y_{2,\epsilon^4 a} = Z_{2,-a^3}$. The
$\wt{A}_{i,a}$ interpolate between the root monomials of
$\U_q(G_2^{(1)})$ and $\U_t(D_4^{(3)})$ as we have the following:

\begin{lem}\label{gi} We have 
$\Pi_q(\wt{A}_{i,a}) = A_{i,\Pi_q(a)}$ for $i\in I, a\in\mathcal{C}$.

We have $\Pi_t(\wt{A}_{2,aq^{-2}}\wt{A}_{2,a}\wt{A}_{2,aq^2}) =
B_{2,(\Pi_t(a))^3}$ for $a\in \mathcal{C}$.

We have $\Pi_t(\wt{A}_{1,a}) = B_{1,-\Pi_t(a)}$ for $a\in \mathcal{C}$.
\end{lem}

\demo The first point is clear.

Let $a' = \Pi_q(a)$. The specialization of
$\wt{A}_{2,aq^{-2}}\wt{A}_{2,a}\wt{A}_{i,aq^2}$ at $q = \epsilon$ is
$$(Y_{2,-a't^{-1}}Y_{2,-a'\epsilon^2 t^{-1}}Y_{2,-a'\epsilon^4
t^{-1}})(Y_{2,-a't}Y_{2,-a\epsilon^2 t}Y_{2,-a'\epsilon^4t}) \times
(Y_{1,a'\epsilon^{-2}}Y_{1,a'}Y_{1,a'\epsilon^2})^{-1}$$
$$ = Z_{2,(a')^3t^3}Z_{2,(a')^3t^{-3}}\times
(Z_{1,a'\epsilon^{-2}}Z_{1,a'}Z_{1,a'\epsilon^2})^{-1} =
B_{2,(a')^3}.$$ The specialization of $\wt{A}_{1,a}$ at $q = \epsilon$
is
$$Y_{1,-a't^{-1}}Y_{1,-a't} \times
 (Y_{2,a'\epsilon^{-2}}Y_{2,a'}Y_{2,a'\epsilon^2})^{-1} =
 Z_{1,-a't^{-1}}Z_{1,-a't}\times Z_{2,(-a')^3}^{-1} = B_{2,-a'}.$$ \qed

Consider the following subalgebras of $\wt{\Yim}_{q,t}$.
\begin{equation*}
\begin{split}
\mathfrak{K}_{1,q,t} = \ZZ[&Y_{1,a}(1 + \wt{A}_{1,aq^3t}^{-1}),
W_{2,a}^{\pm 1}, \beta Y_{2,a}^{\pm 1}, \beta]_{a\in\mathcal{C}}, \\
\mathfrak{K}_{2,q,t} = \ZZ[&Y_{2,a}Y_{2,aq^2}Y_{2,aq^4}(1 + \beta
\wt{A}_{2,aq^5t}^{-1} + \beta
\wt{A}_{2,aq^5t}^{-1}\wt{A}_{2,aq^3t}^{-1}+
\wt{A}_{2,aq^5t}^{-1}\wt{A}_{2,aq^3t}^{-1}\wt{A}_{2,aqt}^{-1}), \\&
\beta Y_{2,a}(1 + \wt{A}_{2,aqt}^{-1}), Y_{1,a}^{\pm 1},
\beta]_{a\in\mathcal{C}}.
\end{split}
\end{equation*}
These are interpolating subalgebras as 

\begin{lem} For $i\in\{1,2\}$, we have $\Pi_q(\mathfrak{K}_{i,q,t}) =
  \mathfrak{K}_{i,q}$ and $\Pi_t(\mathfrak{K}_{i,q,t}) =
  \mathfrak{K}_{i,t}^L$.
\end{lem}

\demo We have
$$\Pi_q(\mathfrak{K}_{2,q,t}) = \ZZ[Y_{2,a}(1 + A_{2,aq}^{-1}),
Y_{1,a}^{\pm 1}]_{a\in q^\ZZ} = \mathfrak{K}_{2,q},$$
$$\Pi_t(\mathfrak{K}_{2,q,t}) =
\ZZ[Y_{2,a}Y_{2,a\epsilon^2}Y_{2,a\epsilon^4}(1
+\wt{A}_{2,-a\epsilon^2t}^{-1}\wt{A}_{2,-a
t}^{-1}\wt{A}_{2,-a\epsilon^{-2} t}^{-1}), Y_{1,a}^{\pm 1}]_{a\in
\epsilon^\ZZ t^\ZZ} = \mathfrak{K}_{2,t}^L,$$ as by Lemma \ref{gi} we
have
$$Y_{2,a}Y_{2,a\epsilon^2}Y_{2,a\epsilon^4}(1
+\wt{A}_{2,-a\epsilon^2t}^{-1}\wt{A}_{2,-a
t}^{-1}\wt{A}_{2,-a\epsilon^{-2} t}^{-1}) = Z_{2,-a^3}(1 +
B_{2,-a^3t^3}^{-1}).$$ Now we have
$$\Pi_q(\mathfrak{K}_{1,q,t}) = \ZZ[Y_{1,a}(1 + A_{1,aq^3}^{-1}),
Y_{2,a}^{\pm 1}]_{a\in q^\ZZ} = \mathfrak{K}_{1,q},$$
$$\Pi_t(\mathfrak{K}_{1,q,t}) = \ZZ[Y_{1,a}(1 + A_{i,-at}^{-1}),
(Y_{2,a}Y_{2,a\epsilon^2}Y_{2,a\epsilon^4})^{\pm 1}]_{a\in
\epsilon^\ZZ t^\ZZ} = \mathfrak{K}_{2,t}^L,$$ as by Lemma \ref{gi} we
have $Y_{1,a}(1 + A_{1,-at}^{-1}) = Z_{1,a} (1 + B_{i, a t}^{-1})$. \qed

As for the case $r = 2$, we define the analogue of $P'\subset P$ in
$\Yim_q$,
\begin{equation*}
\Yim_q' = \ZZ[Y_{1,a}^{\pm 1}, (Y_{2,aq^2}Y_{2,a}Y_{2,aq^{-2}})^{\pm
1}]_{a\in q^\ZZ}.
\end{equation*}

\subsection{Examples}\label{trex}

Now we have to check the existence of interpolating $(q,t)$-characters
in some elementary cases. First, consider the following interpolating
$(q,t)$-character.
$$\xymatrix{&Y_{1,1} \ar[d]^{1,q^3t}&
\\&Y_{1,q^6t^2}^{-1}Y_{2,q^5t}Y_{2,q^3t}Y_{2,qt}\ar[d]^{2,q^6t^2}&
\\&\beta Y_{2,q^7t^3}^{-1}Y_{2,q^3t}Y_{2,qt}\ar[d]^{2,q^4t^2}&
\\&\beta
Y_{1,q^4t^2}Y_{2,q^7t^3}^{-1}Y_{2,q^5t^3}^{-1}Y_{2,qt}\ar[ld]^{2,q^2t^2}\ar[rd]^{1,q^7t^3}&
\\Y_{2,q^7t^3}^{-1}Y_{2,q^5t^3}^{-1}Y_{2,q^3t^3}^{-1}Y_{1,q^4t^2}Y_{1,q^2t^2}\ar[d]^{1,q^5t^3}\ar[rd]^{1,q^7t^3}&&\beta
Y_{1,q^{10}t^4}^{-1}Y_{2,q^9t^3}Y_{2,qt}\ar[ld]^{2,q^2t^2}\ar[d]^{2,q^{10}t^4}
\\Y_{1,q^4t^2}Y_{1,q^8t^4}^{-1}\ar[d]^{1,q^7t^3}&Y_{1,q^2t^2}Y_{1,q^{10}t^4}^{-1}Y_{2,q^9t^3}Y_{2,q^3t^3}^{-1}\ar[ld]^{1,q^5t^3}\ar[rd]^{2,q^{10}t^4}&\beta
Y_{2,q^{11}t^5}^{-1}Y_{2,qt}\ar[d]^{2,q^2t^2}
\\Y_{1,q^8t^4}^{-1}Y_{1,q^{10}t^4}^{-1}Y_{2,q^9t^3}Y_{2,q^7t^3}Y_{2,q^5t^3}\ar[rd]^{2,q^{10}t^4}&&\beta
Y_{1,q^2t^4}Y_{2,q^{11}t^5}^{-1}Y_{2,q^3t^5}^{-1}\ar[ld]^{1,q^5t^5}
\\&\beta
Y_{1,q^8t^4}^{-1}Y_{2,q^{11}t^5}^{-1}Y_{2,q^7t^3}Y_{2,q^5t^3}\ar[d]^{2,q^8t^4}&
\\&\beta
Y_{2,q^{11}t^5}^{-1}Y_{2,q^9t^5}^{-1}Y_{2,q^5t^3}\ar[d]^{2,q^6t^4}&
\\&Y_{2,q^{11}t^5}^{-1}Y_{2,q^9t^5}^{-1}Y_{2,q^7t^5}^{-1}Y_{1,q^6t^4}\ar[d]^{1,q^9t^5}&
\\&Y_{1,q^{12}t^6}^{-1}&}$$
Here we have to check that it is in the $\mathfrak{K}$ as {\em a priori} it is
unclear that
$$\beta
Y_{1,q^{10}t^4}^{-1}Y_{2,q^9t^3}Y_{2,qt}
+Y_{1,q^2t^2}Y_{1,q^{10}t^4}^{-1}Y_{2,q^9t^3}Y_{2,q^3t^3}^{-1}
+ \beta
Y_{2,q^{11}t^5}^{-1}Y_{2,qt}
+ \beta
Y_{1,q^2t^4}Y_{2,q^{11}t^5}^{-1}Y_{2,q^3t^5}^{-1}$$
is in $\mathfrak{K}_{2,q,t}$.
But if we subtract $\beta Y_{2,qt}(1 + A_{2,q^2t^2}^{-1})Y_{2,q^9t^3}(1 +
A_{2,q^{10}t^4}^{-1})Y_{1,q^{10}t^4}^{-1}\in\mathfrak{K}_{2,q,t}$, we get
$$(1 - \beta) Y_{2,q^3t^3}^{-1}Y_{2,q^9t^3}Y_{1,q^{10}t^4}^{-1}Y_{1,q^2t^2}
= (1 - \beta)Y_{1,q^{10}t^4}^{-1}Y_{1,q^2t^2}\in \mathfrak{K}_{2,q,t}.$$
By specializing at $t = 1$, we get the $q$-character of the
15-dimensional fundamental representation of $\U_q(\wh{G_2})$ as
computed in \cite[Appendix]{hadv}. By specializing at $q = \epsilon$, we
get the following.
$$\xymatrix{&Z_{1,1} \ar[d]^{1,t}&
\\&Z_{1,t^2}^{-1} Z_{2,t^3}\ar[ddd]^{2,t^6}&
\\&&
\\&&
\\& Z_{2,t^9}^{-1}Z_{1,\epsilon^4 t^2}Z_{1,\epsilon^2 t^2}\ar[dl]^{1,-\epsilon^5 t^3}\ar[rd]^{1,-\epsilon t^3}&
\\Z_{1,\epsilon^4 t^2}Z_{1,\epsilon^2 t^4}^{-1}\ar[dr]^{1,-\epsilon t^3}&&Z_{1,\epsilon^2 t^2}Z_{1,\epsilon^4 t^4}^{-1}\ar[ld]^{1,-\epsilon^5 t^3}
\\&Z_{1,\epsilon^2 t^4}^{-1}Z_{1,\epsilon^4 t^4}^{-1}Z_{2,t^9}\ar[ddd]^{2,t^{12}}&
\\&&
\\&&
\\&Z_{2,t^{15}}^{-1} Z_{1,t^4}\ar[d]^{1,t^5}&
\\&Z_{1,t^6}^{-1}&}$$
This is the twisted $t$-character of the $8$-dimensional fundamental
representation of $\U_t(D_4^{(3)})$ as computed in \cite[Section 11.2]{hgen}.

Now we have to consider the case of the monomial
$Y_{2,1}Y_{2,q^2}Y_{2,q^4} = W_{2,q^2}$. The dimension of the
corresponding KR module of $\U_q(G_2^{(1)})$ is
$133$ (this can be obtained, for example, from the $T$-system proved
in \cite{hcr}: let $T_k^{(i)}$ be the dimension of a
KR module of highest weight $k\omega_i$. Then for the
fundamental representations we have $T_1^{(1)} = 15$, $T_1^{(2)} = 7$,
so $T_2^{(2)} = (T_1^{(2)})^2 - T_1^{(1)} = 34$ and $T_3^{(2)} =
(T_1^{(2)})^{-1}((T_2^{(2)})^2 - (T_1^{(1)})^2) = 133$).

There is also an interpolating $(q,t)$-character in this case. We do
not list all $133$ monomials, but we list the $29$ monomials (with multiplicity) 
which do not have $\beta$ in their coefficient:

$2_1 2_{q^2} 2_{q^4} ;
2_{q^2t^2}^{-1}2_{q^4t^2}^{-1}2_{q^6t^2}^{-1}1_{q^5t}1_{q^3t}1_{qt};
1_{q^3t}1_{qt}1_{q^{11}t^3}^{-1}2_{q^2t^2}^{-1}2_{q^4t^2}^{-1}2_{q^8t^2}2_{q^{10}t^2};
1_{q^5t}1_{q^9t^3}^{-1}1_{qt}2_{q^2t^2}^{-1}2_{q^8t^2};
\\1_{q^5t}1_{q^3t}1_{q^7t^3}^{-1};
1_{q^9t^3}^{-1}1_{qt}1_{q^{11}t^3}^{-1}2_{q^2t^2}^{-1}2_{q^6t^2}2_{q^8t^2}^22_{q^{10}t^2};
1_{q^5t}1_{q^9t^3}^{-1}1_{q^7t^3}^{-1}2_{q^4t^2}2_{q^6t^2}2_{q^8t^2};
\\1_{q^{11}t^3}^{-1}1_{q^3t}1_{q^7t^3}^{-1}2_{q^{10}t^2}2_{q^8t^2}2_{q^6t^2};
1_{qt}1_{q^7t^3}2_{q^2t^2}^{-1}2_{q^8t^4}^{-1}2_{q^8t^2}2_{q^{10}t^4}^{-1}2_{q^{12}t^4}^{-1};
1_{q^5t}1_{q^5t^3}2_{q^6t^4}^{-1}2_{q^8t^4}^{-1}2_{q^{10}t^4}^{-1};
\\1_{q^3t}1_{q^9t^3}2_{q^{12}t^4}^{-1}2_{q^{10}t^4}^{-1}2_{q^8t^4}^{-1};
1_{q^9t^3}^{-1}1_{q^7t^3}^{-1}1_{q^{11}t^3}^{-1}2_{q^4t^2}2_{q^6t^2}^22_{q^8t^2}^22_{q^{10}t^2};
1_{qt}1_{q^{13}t^5}^{-1}2_{q^2t^2}^{-1}2_{q^8t^2};
1_{q^5t}1_{q^{11}t^5}^{-1};
\\1_{q^3t}1_{q^{15}t^5}^{-1}2_{q^{14}t^4}2_{q^8t^4}^{-1};
2\times 2_{q^4t^2}2_{q^6t^2}2_{q^{10}t^4}^{-1}2_{q^{12}t^4}^{-1};
1_{q^7t^3}^{-1}1_{q^{13}t^5}^{-1}2_{q^4t^2}2_{q^6t^2}2_{q^8t^2};
\\1_{q^{11}t^3}^{-1}1_{q^{11}t^5}^{-1}2_{q^{10}t^2}2_{q^8t^2}2_{q^6t^2};
1_{q^9t^3}^{-1}1_{q^{15}t^5}^{-1}2_{q^{14}t^4}2_{q^6t^4}2_{q^4t^4};
2_{q^6t^4}^{-1}2_{q^8t^4}^{-2}2_{q^{10}t^4}^{-2}2_{q^{12}t^4}^{-1}1_{q^5t^3}1_{q^7t^3}1_{q^9t^3};
\\2_{q^6t^4}^{-1}2_{q^8t^4}^{-1}2_{q^{10}t^4}^{-1}2_{q^{14}t^4}1_{q^5t^3}1_{q^7t^3}1_{q^{15}t^5}^{-1};
2_{q^6t^4}^{-1}2_{q^8t^4}^{-1}2_{q^{10}t^4}^{-1}1_{q^5t^3}1_{q^{13}t^5}^{-1}1_{q^9t^3};
\\2_{q^8t^4}^{-1}2_{q^{10}t^4}^{-1}2_{q^{12}t^4}^{-1}1_{q^{11}t^5}^{-1}1_{q^7t^3}1_{q^9t^3};
2_{q^6t^4}^{-1}2_{q^{12}t^4}2_{q^{14}t^4}1_{q^5t^3}1_{q^{13}t^5}^{-1}1_{q^{15}t^5}^{-1};
1_{q^{11}t^5}^{-1}1_{q^{13}t^5}^{-1}1_{q^9t^3};
\\2_{q^8t^4}^{-1}2_{q^{14}t^4}1_{q^{11}t^5}^{-1}1_{q^7t^3}1_{q^{15}t^5}^{-1};
1_{q^{11}t^5}^{-1}1_{q^{13}t^5}^{-1}1_{q^{15}t^5}^{-1}2_{q^{14}t^4}2_{q^{12}t^4}2_{q^{10}t^4};
2_{q^{16}t^6}^{-1}2_{q^{14}t^6}^{-1}2_{q^{12}t^6}^{-1}.$
  
As the other terms disappear when we specialize at $q = \epsilon$, we
can compute the specialization from the above terms which is given
in the figure bellow. We get
the twisted $t$-character of the $29$-dimensional fundamental
representation of $\U_t(D_4^{(3)})$ as computed in \cite[Section 11.2]{hgen}.

\begin{figure}[htbp]
\begin{center}
\leavevmode
\begin{equation*}
\divide\dgARROWLENGTH by 3
\dgARROWPARTS=6
\begin{diagram}
  \node[3]{2_{-1}} \arrow{s,r}{2,-t^3}
\\
  \node[3]{2_{-t^6}^{-1} 1_{\epsilon t} 1_{- t} 1_{\epsilon^5 t}}
  \arrow{sw,t,3}{1,\epsilon t^2} \arrow{s,l,3}{1,- t^2}
  \arrow{se,t,3}{1,\epsilon^5 t^2}
\\
  \node[2]{1_{\epsilon t^3}^{-1} 1_{- t} 1_{\epsilon^5 t}}
  \arrow{s,l}{1,- t^2} \arrow{se,b,1}{1,\epsilon^5 t^2}
  \node{1_{\epsilon t} 1_{- t^3}^{-1} 1_{\epsilon^5 t}}
  \arrow{sw,b,1}{1,\epsilon t^2}\arrow{se,b,1}{1,\epsilon^5 t^2}
  \node{1_{\epsilon t} 1_{- t} 1_{\epsilon^5 t^3}^{-1}}
  \arrow{sw,b,1}{1, \epsilon t^2} \arrow{s,r}{1, - t^2}
\\
  \node[2]{1_{\epsilon t^3}^{-1} 2_{-t^6} 1_{- t^3}^{-1} 1_{\epsilon^5 t}}
  \arrow{sw,b,1}{2,-t^9} \arrow{se,b,1}{1,\epsilon^5 t^2}
  \node{1_{\epsilon t^3}^{-1} 2_{-t^6} 1_{- t} 1_{\epsilon^5 t^3}^{-1}}
  \arrow{sw,b,1}{2,-t^9}
  \arrow{s,r}{1,- t^2}
  \node{1_{\epsilon t} 2_{- t^6} 1_{- t^3}^{-1} 1_{\epsilon^5 t^3}^{-1}}
  \arrow{sw,b,1}{1, \epsilon t^2}
  \arrow{s,r}{2, -t^9}
\\
  \node{2_{-t^{12}}^{-1}1_{\epsilon^5 t}1_{\epsilon^5 t^3}}
  \arrow{s,l}{1,\epsilon^5 t^4}
  \node{1_{- t} 1_{- t^3} 2_{-t^{12}}^{-1}}
  \arrow{s,r,1}{1, - t^4}
  \node{2_{-t^6}^2 1_{\epsilon t^3}^{-1} 1_{- t^3}^{-1} 1_{\epsilon^5 t^3}^{-1}}
  \arrow{s,r}{2, -t^9}
  \node{2_{-t^{12}}^{-1} 1_{\epsilon t} 1_{\epsilon t^3}}
  \arrow{s,r}{1, \epsilon t^4}
\\
  \node{1_{\epsilon^5 t}1_{\epsilon^5 t^5}^{-1}}
  \arrow{s,l}{1, \epsilon^5 t^2}
  \node{1_{- t}1_{- t^3}^{-1}}
  \arrow{s,l}{1, - t^2}
  \node{2\times 2_{-t^6}2_{-t^{12}}^{-1}}
  \arrow{s,r}{2, -t^9}
  \node{1_{\epsilon t} 1_{\epsilon t^5}^{-1}}
  \arrow{s,r}{1, \epsilon t^2}
\\
  \node{2_{-t^6}1_{\epsilon^5 t^3}^{-1}1_{\epsilon^5 t^5}^{-1}}
  \arrow{se,b,1}{2, -t^9}
  \node{1_{- t^3}^{-1}1_{- t^5}^{-1} 2_{-t^6}}
  \arrow{se,b,1}{2, -t^9}
  \node{1_{\epsilon t^3}2_{-t^{12}}^{-2}1_{- t^3}1_{\epsilon^5 t^3}}
  \arrow{sw,b,1}{1, \epsilon^5 t^4}\arrow{s,r,2}{1, - t^4}
  \arrow{se,t,3}{1, \epsilon t^4}
  \node{2_{-t^6}1_{\epsilon t^3}^{-1}1_{\epsilon t^5}^{-1}}
  \arrow{s,r}{2, -t^9}
\\
  \node[2]{1_{- t^3}  2_{-t^{12}}^{-1} 1_{\epsilon t^3} 1_{\epsilon^5 t^5}^{-1}}
  \arrow{s,l}{1, - t^4}
  \arrow{se,b,1}{1,\epsilon t^4}
  \node{1_{- t^5}^{-1} 2_{-t^{12}}^{-1} 1_{\epsilon t^3} 1_{\epsilon^5 t^3}}
  \arrow{sw,b,1}{1,\epsilon^5 t^4}\arrow{se,b,1}{1, \epsilon t^4}
  \node{1_{\epsilon^5 t^3}  2_{-t^{12}}^{-1}1_{\epsilon t^5}^{-1}1_{- t^3}}
  \arrow{sw,b,1}{1, \epsilon^5 t^4} \arrow{s,r}{1,- t^4}
\\
  \node[2]{1_{- t^5}^{-1} 1_{\epsilon t^3} 1_{\epsilon^5 t^5}^{-1}}
  \arrow{se,b,2}{1, \epsilon t^4}
  \node{1_{- t^3} 1_{\epsilon t^5}^{-1} 1_{\epsilon^5 t^5}^{-1}}
  \arrow{s,r}{1, - t^4}
  \node{1_{\epsilon t^5}^{-1} 1_{- t^5}^{-1} 1_{\epsilon^5 t^3}}
  \arrow{sw,b,2}{1, \epsilon^5 t^4}
\\
  \node[3]{1_{\epsilon t^5}^{-1}2_{-t^{12}}1_{- t^5}^{-1}1_{\epsilon^5 t^5}^{-1}}
  \arrow{s,r}{2, -t^{15}}
\\
  \node[3]{2_{-t^{18}}^{-1}}
\end{diagram}
\end{equation*}
\end{center}
\end{figure}

\newpage

\subsection{Conclusion}

With the existence of the two elementary interpolating
$(q,t)$-characters in the last subsection, we can conclude the proof of
the two main results of this section. We define $\mathfrak{K}$ as for
the case $r = 2$ and we have the following:

\begin{thm}\label{virtg} For all dominant $m$ there is a unique
  $F(m)\in\mathfrak{K}$ such that $m$ is the unique dominant
  monomial of $F(m)$.\end{thm}

As in the double-laced case, we have the notion of Langlands dual
representation and interpolating $(q,t)$-character in $\mathfrak{K}$
with highest monomial in $\wt{\Yim}_{q,t}\setminus \beta
  \wt{\Yim}_{q,t}$. We get the following consequence of Theorem \ref{virtg}:

\begin{thm}\label{dualg} An irreducible tensor product of KR modules over $\U_q(G_2^{(1)})$
  of highest monomial in $\Yim_q'$ has a Langlands dual
  representation. Moreover, the Langlands dual representation of a KR
  module over $\U_q(G_2^{(1)})$ is a KR module over
  $\U_t(D_4^{(3)})$.
\end{thm}

\section{From twisted to untwisted types}\label{reverse}

In this section we describe the Langlands duality in the opposite
direction, from a twisted quantum algebra $\U_t({}^L\hat{\Glie})$ to
an untwisted quantum affine algebra $\U_q(\hat{\Glie})$. We prove the
existence of interpolating $(t,q)$-characters and we prove the duality
for irreducible tensor product of KR modules (for this duality we have to use a
slightly generalized definition of KR modules over
twisted quantum affine algebras).

\subsection{Double-laced cases}

We use the notation of Section \ref{double}, in particular, for
$\phi$, $\epsilon$.  Note that $q_i = q^{r_i}$ and not $q^{r_i^\vee}$.
We need the function $\alpha^L(t,q)$ such that $\alpha^L(t,\epsilon) =
1$ and $\alpha^L(1, q) = 0$ defined by $\alpha^L(t,q) = 1 -
\alpha(q,t)$.  Consider the ring
$$\Yim_{t,q}^L = \ZZ[X_{i,a}^{\pm 1}, \alpha^L z_{i,a}^{\pm 1},
  \alpha^L]_{i\in I, a\in\mathcal{C}}\subset \ZZ[z_{i,a}^{\pm
    1},\alpha^L]_{i\in I, a\in\mathcal{C}},$$
\begin{equation*}
\begin{split}
\text{ where }X_{i,a} = \begin{cases} z_{i,a}&\text{ if $i\in I_2^\vee = I_1$,}
                     \\z_{i,aq^{-1}}z_{i,aq}&\text{ if $i\in I_1^\vee = I_2$.}
\end{cases}
\end{split}
\end{equation*}
We then have surjective specialization maps, respectively, at $q =
\epsilon$ and $t = 1$,
$$\Pi_t^L : \Yim_{t,q}^L\rightarrow
  \ZZ[Z_{i,a^{r_i^{\vee}}}^{\pm 1}]_{i\in I, a\in \epsilon^\ZZ t^\ZZ}
  = \Yim_t^L,$$
$$\Pi_q^L : \Yim_{t,q}^L\rightarrow \ZZ[Y_{i,a}^{\pm 1}]_{i\in I, a\in
  q^\ZZ} = \Yim_q,$$
where for $a\in \mathcal{C}$, $i\in I$, we assign
$$X_{i,a} \mapsto Y_{i,a}\text{   and   }z_{i,a} \mapsto Z_{i,(-1)^{1 + \phi(i)}(a)^{r_i^\vee}}.$$
Note that for $i\in I_2^\vee, a\in \epsilon^\ZZ t^\ZZ$, $\Pi_t^L(z_{i,a}) = \Pi_t^L(z_{i,-a}) = Z_{i,a^2}$.
We have the ideals
$$\text{Ker}(\Pi_t^L) = \langle (\alpha^L - 1), (X_{i,aq} - X_{i,a\epsilon}),
\alpha^L (z_{i,aq} - z_{i,a\epsilon}), (z_{j,a} - z_{j,-a})\rangle_{i\in I, j\in I_2^\vee,a\in\mathcal{C}},$$
$$\text{Ker}(\Pi_q^L) = \langle \alpha^L, (X_{i,at} - X_{i,a})
  \rangle_{i\in I, a\in\mathcal{C}}.$$ 
  The ideal $\text{Ker}(\Pi_t^L)\cap \text{Ker}(\Pi_q^L)$ is generated
  by the elements
$$\alpha^L (\alpha^L - 1)\text{ , }\alpha^L (z_{i,aq} -
  z_{i,a\epsilon})\text{ , }\alpha^L (z_{j,a} - z_{j,-a}),$$
$$(\alpha^L - 1) (X_{i,at} - X_{i,a})\text{ , }(X_{i,a} -
  X_{i,at})(X_{k,bq} - X_{k,b\epsilon})\text{ , }(X_{i,a} -
  X_{i,at})(z_{j,b} - z_{j,-b}),$$
for $i, k\in I$, $j\in I_2^\vee$, $a, b\in\mathcal{C}$.
We will work in the ring
$$\wt{\Yim}_{t,q}^L = \Yim_{t,q}^L/(\text{Ker}(\Pi_t^L)\cap
\text{Ker}(\Pi_q^L)).$$
We use the notion of monomial, dominant monomial as above. 

\begin{defi} For $i\in I$ and $a\in \mathcal{C}$ we define
$$\wt{B}_{i,a} = z_{i,a(q_it)^{-1}}z_{i,aq_it} \times \prod_{j\in
I,C_{j,i} = - 1}z_{j,a}^{-1}\times \prod_{j\in I,C_{j,i} = -
2}z_{j,aq^{-1}}^{-1} z_{j,aq}^{-1}.$$
\end{defi}

\begin{lem} We have $\Pi_t^L(\wt{B}_{i,a}) =
  B_{i,(\Pi_t^L(a))^{r_i^\vee}(-1)^{\phi(i)}}$ for $i\in I$,
  $a\in\mathcal{C}$.

We have $\Pi_q^L(\wt{B}_{i,aq^{-1}}\wt{B}_{i,aq}) =
A_{i,\Pi_q^L(a)}$ for $i\in I_1^\vee$, $a\in\mathcal{C}$.

We have $\Pi_q^L(\wt{B}_{i,a}) =
A_{i,\Pi_q^L(a)}$ for $i\in I_2^\vee$, $a\in\mathcal{C}$.
\end{lem}

\demo Let $a' = \Pi_t^L(a)$. For $i\in I_1^\vee$, we have
$$\Pi_t^L(\wt{B}_{i,a}) = z_{i,-a' t^{-1}}z_{i,-a't} \times \prod_{j\in
I,C_{j,i} = - 1}z_{j,a'}^{-1}$$
$$= Z_{i,(-1)^{\phi(i)}a' t^{-1}}Z_{i,(-1)^{\phi(i)}a' t} \times \prod_{j\in
I_1^{\vee},C_{j,i} = - 1}Z_{j,(-1)^{\phi(i)}a'}^{-1}\times \prod_{j\in
I_2^{\vee},C_{j,i} = - 1}Z_{j,(a')^2}^{-1},$$
which is equal to $B_{i,a'(-1)^{\phi(i)}}$. Indeed if there is $j\in
I_2^\vee = I_1$ such that $C_{j,i} = -1$, we have $\phi(j) = 1$
and so $Z_{j,(a')^2} = z_{j,a'}$.

For $i\in I_2^\vee$, we have
$$\Pi_t^L(\wt{B}_{i,a}) = z_{i,-\epsilon a' t^{-1}}z_{i,\epsilon a't}
\times \prod_{j\in
I,C_{j,i} = - 1}z_{j,a'}^{-1}\times \prod_{j\in
I,C_{j,i} = - 2}z_{j,-\epsilon a'}^{-1}z_{j,\epsilon a'}^{-1}$$
$$= Z_{i,(-1)^{\phi(i)}(a')^2 t^{-2}}Z_{i,(-1)^{\phi(i)}(a')^2 t^2}
\times \prod_{j\in
I,C_{j,i} = - 1}Z_{j,(-1)^{\phi(i)}(a')^2}^{-1}\times \prod_{j\in
I,C_{j,i} = - 2}Z_{j,\epsilon a'}^{-1}Z_{j,- \epsilon a'}^{-1},$$ 
which is equal to $B_{i,(-1)^{\phi(i)} (a')^2}$.
Indeed if there is $j\in I$ such that $C_{j,i} = -2$, we have $\phi(i)
= 1$ and so $(\epsilon a')^2 = - (a')^2 = (-1)^{\phi(i)} (a')^2$.

Let $a'' = \Pi_q^L(a)$. For $i\in I_1^\vee$, we have
$$\Pi_q^L((\wt{B}_{i,aq^{-1}}\wt{B}_{i,aq})) 
= (z_{i,a''q^{-3}}z_{i,a''q^{-1}})(z_{i,a''q}z_{i,a''q^3}) 
\times \prod_{j\in I,C_{j,i} = - 1}z_{j,a''q^{-1}}^{-1}z_{j,a''q}^{-1}$$
$$= Y_{i,a''q^{-2}} Y_{i,a''q^2} 
\times \prod_{j\in I_1^{\vee},C_{j,i} = - 1}Y_{j,a''}^{-1}
\times \prod_{j\in I_2^{\vee},C_{j,i} = -
  1}Y_{j,a''q^{-1}}^{-1}Y_{j,a''q}^{-1} = A_{i,a''}.$$
For $i\in I_2^\vee$, we have
$$\Pi_q^L(\wt{B}_{i,a}) = z_{i,a'' q^{-1}}z_{i,a''q} \times \prod_{j\in
I,C_{j,i} = - 1}z_{j,a''}^{-1}\times \prod_{j\in I,C_{j,i} = -
2}(z_{j,a''q^{-1}} z_{j,a''q})^{-1}$$
$$= Y_{i,a'' q^{-1}}Y_{i,a''q} \times \prod_{j\in
I,C_{j,i} < 0}Y_{j,a''}^{-1} = A_{i,a''}.$$
\qed

For $i\in I_1^\vee$, consider the subalgebra $\mathfrak{K}_{i, t,q}^L$
of $\Yim_{t,q}^L$ equal to
$$\ZZ[X_{i,a}(1 + \alpha^L
\wt{B}_{i,aq^3t}^{-1} + \alpha^L \wt{B}_{i,aqt}^{-1} 
+ \wt{B}_{i,aq^3t}^{-1}\wt{B}_{i,aqt}^{-1}),
\alpha^L z_{i,a}(1 + \wt{B}_{i,aq^2t}^{-1}), X_{j,a}^{\pm 1}, \alpha^L
z_{j,a}^{\pm 1}, \alpha^L]_{a\in\mathcal{C}, j\neq i},$$ 
and for $i\in I_2^L$,
$$\mathfrak{K}_{i,t,q}^L = \ZZ[z_{i,a}(1 + \wt{B}_{i,aqt}^{-1}),
  X_{j,a}^{\pm 1}, \alpha^L z_{j,a}^{\pm 1}, \alpha^L]_{a\in
  \mathcal{C}, j\neq i}.$$
Then we have the following.

\begin{lem} We have 
$\Pi_t^L(\mathfrak{K}_{i,t,q}^L)= \mathfrak{K}_{i,t}^L$ and 
$\Pi_q^L(\mathfrak{K}_{i,t,q}^L) = \mathfrak{K}_{i,q}$ for $i\in I$.
\end{lem}

\demo For $i\in I_1^\vee$, $\Pi_t^L(\mathfrak{K}_{i, t,q}^L)$ is equal to
$$\ZZ[z_{i,-\epsilon a}z_{i,\epsilon a}(1 + 
B_{i,-\epsilon at (-1)^{\phi(i)}}^{-1})( 1+ B_{i,\epsilon at
  (-1)^{\phi(i)}}^{-1}),
z_{i,a}(1 + B_{i,a t(-1)^{\phi(i) + 1}}^{-1}), z_{j,a}^{\pm 1}]_{a\in
  t^{\ZZ}\epsilon^\ZZ, j\neq i}$$
$$= \ZZ[Z_{i,a(-1)^{\phi(i) + 1}}(1 + B_{i,a (-1)^{\phi(i) +
      1}t}^{-1}), Z_{j,a^{r_j^\vee}}^{\pm 1}]_{a\in
  t^{\ZZ}\epsilon^\ZZ, j\neq i} = \mathfrak{K}_{i,t}^L.$$
We also have
$$\Pi_q^L(\mathfrak{K}_{i, t,q}^L) = \ZZ[\Pi_q^L(X_{i,a}(1 +
  \wt{B}_{i,aq^3t}^{-1}\wt{B}_{i,aqt}^{-1})),
\Pi_q^L(X_{j,a}^{\pm 1})]_{a\in\mathcal{C}, j\neq i}$$
$$= \ZZ[Y_{i,a}(1 + A_{i,aq^2}^{-1}), Y_{j,a}^{\pm 1}]_{a\in q^{\ZZ},
  j\neq i} = \mathfrak{K}_{i,q}.$$
Now for $i\in I_2^\vee$, we have
$$\Pi_t^L(\mathfrak{K}_{i,t,q}^L) = \ZZ[z_{i,a}(1 + B_{i,a^2
  t^2(-1)^{\phi(i) + 1}}^{-1}),
  z_{j,a}^{\pm 1}]_{a\in t^\ZZ\epsilon^\ZZ, j\neq i}$$
  $$= \ZZ[Z_{i,a^2(-1)^{\phi(i) + 1}}(1 + B_{i,a^2 t^2(-1)^{\phi(i) +
  1}}^{-1}),
  Z_{j,a^{r_j^\vee}}^{\pm 1}]_{a\in t^\ZZ\epsilon^\ZZ, j\neq i} =
  \mathfrak{K}_{i,t}^L,$$
$$\Pi_q^L(\mathfrak{K}_{i,t,q}^L) = \ZZ[Y_{i,a}(1 + A_{i,aq}^{-1}),
  Y_{j,a}^{\pm 1}]_{a\in q^\ZZ, j\neq i} = \mathfrak{K}_{i,q}.$$ 
\qed

We set
\begin{equation*}
(\Yim_t^L)' = \ZZ[Z_{i,a}^{\pm 1}]_{i \in I_2^\vee,a\in (\epsilon^\ZZ
    t^\ZZ)^2}\otimes
\ZZ[(Z_{i,a}Z_{i,-a})^{\pm 1}]_{i \in I_1^\vee,a\in \epsilon^\ZZ
    t^\ZZ}\subset \Yim_t^L.
\end{equation*}
and we define $\mathfrak{K}^L\subset \tilde{\Yim}_{t,q}^L$ as above.

As in the previous sections, we check the existence of various
elements in $\mathfrak{K}^L$ that we call interpolating
$(t,q)$-characters.

First, we suppose that $\U_t({}^L\wh{\Glie})$ is of type $A_3^{(2)}$,
and so that $\U_q(\wh{\Glie})$ is of type $C_2^{(1)}$, with $r_1 =
r_2^\vee = 2$ and $r_2 = r_1^\vee = 1$.  We have $\phi(1) = 0$,
$\phi(2) = 1$. We have the following interpolating $(q,t)$-character.
$$\xymatrix{z_{2,1} \ar[d]^{2,qt} &&Z_{2,1} \ar[d]^{2,t^2} & &  Y_{2,1}\ar[d]^{2,q}
\\ z_{2,q^2t^2}^{-1} z_{1,t}z_{1,tq^2}\ar[d]^{1,q^4 t^2}\ar[dr]^{1,q^2 t^2}&&Z_{2,t^4}^{-1} Z_{1,-t}Z_{1,t}\ar[d]^{1, t^2}\ar[dr]^{1,- t^2}& & Y_{2,q^2}^{-1}Y_{1,q}\ar[dd]^{1,q^2}
\\ \alpha^L z_{1,t^3q^6}^{-1}z_{1,t}z_{2,t^2q^4}z_{2,t^2q^2}^{-1}\ar[d]^{1,q^2t^2}
& \alpha^L z_{1,tq^2}z_{1, t^3q^4}^{-1} \ar[dl]^{1,q^4t^2}&Z_{1,t^3}^{-1}Z_{1,-t} \ar[d]^{1,t^2}
&  Z_{1,t}Z_{1, -t^3}^{-1} \ar[dl]^{1,-t^2}&
\\ z_{1,t^3q^4}^{-1}z_{1,t^3q^6}^{-1} z_{2,t^2q^4}\ar[d]^{2,t^3q^5}&&Z_{1,-t^3}^{-1}Z_{1,t^3}^{-1} Z_{2,t^4}\ar[d]^{2,t^6}&&Y_{1,q^5}^{-1}Y_{2,q^4}\ar[d]^{2,q^5} 
\\z_{2,t^4q^6}^{-1}&&Z_{2,t^8}^{-1}&&Y_{2,q^6}^{-1}}$$
$\Pi_t^L$ gives the twisted $t$-character of a $6$-dimensional fundamental
representation of $\U_t(A_3^{(2)})$ and $\Pi_q^L$ the $q$-character
of a $4$-dimensional fundamental representation of $\U_q(C_2^{(1)})$. 

We also have the following interpolating $(t,q)$-character
$$\xymatrix{&z_{1,t}z_{1,tq^2}\ar[dl]^{1,t^2q^4}\ar[dr]^{1,t^2q^2}&
\\ \alpha^Lz_{2,t^2q^4}z_{1,t}z_{1,t^3q^6}^{-1}\ar[d]^{2,t^3q^5}\ar[dr]^{1,t^2q^2}&&\alpha^L z_{2,t^2q^2}z_{1,tq^2}z_{1,t^3q^4}^{-1}\ar[d]^{2,t^3q^3}\ar[dl]^{1,t^2q^4}
\\ \alpha^L z_{2,t^4q^6}^{-1}z_{1,t}z_{1,t^3q^4} \ar[d]^{1,t^4q^6}& z_{2,t^2q^4}z_{2,t^2q^2}z_{1,t^3q^6}^{-1}z_{1,t^3q^4}^{-1}\ar[d]^{2,t^3q^5}
&\alpha^L z_{2,t^4q^4}^{-1}z_{1,tq^2}z_{1,t^3q^2}\ar[d]^{1,t^4q^4}
\\ \alpha^L z_{1,t}z_{1,t^5q^8}^{-1} \ar[d]^{1,t^2q^2}& ( 1 + \alpha^L) z_{2,t^4q^6}^{-1}z_{2,t^2q^2}\ar[d]^{2,t^3q^3}&\alpha^L z_{1,tq^2}z_{1,t^5q^6}^{-1}\ar[d]^{1,t^2q^4}
\\ \alpha^L z_{2,t^2q^2}z_{1,t^3q^4}^{-1}z_{1,t^5q^8}^{-1}\ar[d]^{2,t^3q^3} &z_{2,t^4q^6}^{-1}z_{2,t^4q^4}^{-1}z_{1,t^3q^2}z_{1,t^3q^4}\ar[dl]^{1,t^4q^6}\ar[dr]^{1,t^4q^4} &\alpha^L z_{2,t^2q^4}z_{1,t^5q^6}^{-1}z_{1,t^3q^6}^{-1} \ar[d]^{2,t^3q^5} 
\\\alpha^L z_{1,t^3q^2}z_{1,t^5q^8}^{-1}z_{2,t^4q^4}^{-1}\ar[dr]^{1,t^4q^4} & & \alpha^L z_{1,t^3q^4}z_{1,t^5q^6}^{-1}z_{2,t^4q^6}^{-1}\ar[dl]^{1,t^4q^6}
\\ &z_{1,t^5q^8}^{-1}z_{1,t^5q^6}^{-1}&&}$$
It is easy to check that it is in the $\mathfrak{K}^L$, for example
$$z_{2,t^2q^4}z_{2,t^2q^2}z_{1,t^3q^6}^{-1}z_{1,t^3q^4}^{-1}
+(1 + \alpha^L) z_{2,t^4q^6}^{-1}z_{2,t^2q^2}
+z_{2,t^4q^6}^{-1}z_{2,t^4q^4}^{-1}z_{1,t^3q^2}z_{1,t^3q^4}$$
$$= 
(1 - \alpha^L) 
+ z_{2,t^2q^4}z_{2,t^2q^2}z_{1,t^3q^6}^{-1}z_{1,t^3q^4}^{-1}(1 + \wt{B}_{2,t^3q^5})(1 + \wt{B}_{t^3q^3})
\in \mathfrak{K}_{2,q,t}^L.$$
Note that the coefficients $\alpha^L$ are imposed by the condition
that the interpolating $(t,q)$-character is in $\mathfrak{K}^L$, in
particular the coefficient $(1 + \alpha^L)$ of
$z_{2,t^4q^6}^{-1}z_{2,t^2q^2}$.

$\Pi_t^L$ gives the twisted $t$-character of a tensor product of two
$4$-dimensional fundamental
representation of $\U_t(A_3^{(2)})$ and $\Pi_q^L$ the $q$-character
of a $5$-dimensional fundamental representation of $\U_q(C_2^{(1)})$. 
$$\xymatrix{&Z_{1,-t}Z_{1,t}\ar[dl]^{1,t^2}\ar[dr]^{1,-t^2}&&Y_{1,q}\ar[dd]^{1,q^3}
\\ Z_{2,t^4}Z_{1,-t}Z_{1,t^3}^{-1}\ar[d]^{2,t^6}\ar[dr]^{1,-t^2}&& Z_{2,t^4}Z_{1,t}Z_{1,-t^3}^{-1}\ar[d]^{2,t^6}\ar[dl]^{1,t^2}&
\\  Z_{2,t^8}^{-1}Z_{1,-t}Z_{1,-t^3} \ar[d]^{1,-t^4}\ar[dr]^{1,-t^2}& Z_{2,t^4}^2Z_{1,t^3}^{-1}Z_{1,-t^3}^{-1}\ar[d]^{2,t^6}
& Z_{2,t^8}^{-1}Z_{1,t}Z_{1,t^3}\ar[dl]^{1,t^2}\ar[d]^{1,t^4}&Y_{1,q^5}^{-1}Y_{2,q^2}Y_{2,q^4}\ar[d]^{2,q^5}
\\  Z_{1,-t}Z_{1,-t^5}^{-1} \ar[d]^{1,-t^2}&2\times Z_{2,t^8}^{-1}Z_{2,t^4}\ar[dl]^{1,-t^4}\ar[d]^{2,t^6}\ar[dr]^{1,t^4}& Z_{1,t}Z_{1,t^5}^{-1}\ar[d]^{1,t^2}&Y_{2,q^2}Y_{2,q^6}^{-1}\ar[d]^{2,q^3}
\\  Z_{2,t^4}Z_{1,-t^3}^{-1}Z_{1,-t^5}^{-1}\ar[d]^{2,t^6} &Z_{2,t^8}^{-2}Z_{1,t^3}Z_{1,-t^3}\ar[dl]^{1,-t^4}\ar[dr]^{1,t^4} & Z_{2,t^4}Z_{1,t^5}^{-1}Z_{1,t^3}^{-1}\ar[d]^{2,t^6}&Y_{1,q^3}Y_{2,q^4}^{-1}Y_{2,q^6}^{-1}\ar[dd]^{1,q^5}  
\\ Z_{1,t^3}Z_{1,-t^5}^{-1}Z_{2,t^8}^{-1}\ar[dr]^{1,t^4} & &  Z_{1,-t^3}Z_{1,t^5}^{-1}Z_{2,t^8}^{-1}\ar[dl]^{1,-t^4}&
\\ &Z_{1,t^5}^{-1}Z_{1,-t^5}^{-1}&&Y_{1,q^7}^{-1}}$$
The multiplicity $2$ of $Z_{2,t^4}Z_{2,t^8}^{-1}$ in the image by $\Pi_t^L$ is ramified into $1 + \alpha^L$ in the interpolating $(t,q)$-character. That is why we get just a multiplicity $1$ for $Y_{2,q^2}Y_{2,q^6}^{-1}$. Note that in particular the interpolating $(t,q)$-character can not be factorized.

Next, consider $C_{i,j} = C_{j,i} = -1$ with $r_1^\vee = r_2^\vee = 1$. We choose $\phi(1) = 0$,
$\phi(2) = 1$, and we have the following interpolating $(t,q)$-character. 
$$
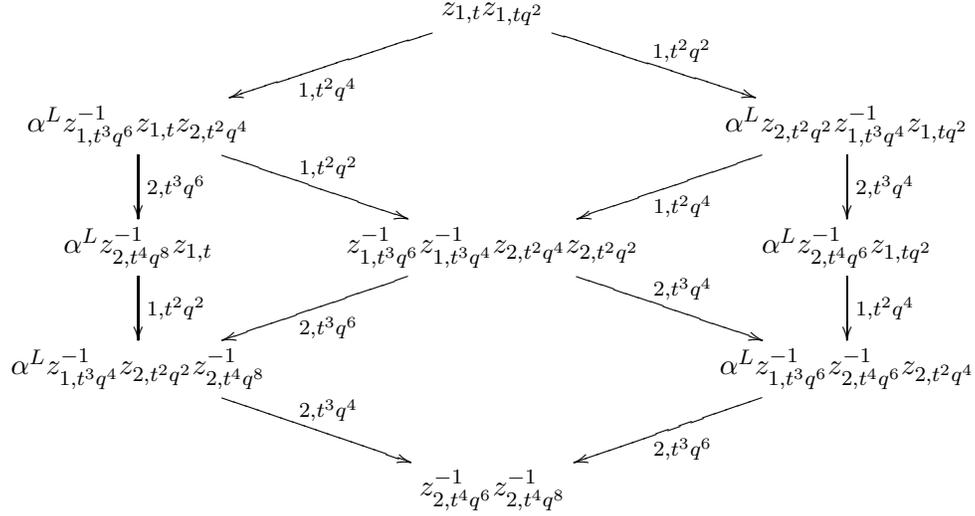
\xymatrix{&z_{1,t}z_{1,tq^2}\ar[dl]^{1,t^2q^4}\ar[dr]^{1,t^2q^2}&
\\\alpha^L z_{1,t^3q^6}^{-1}z_{1,t}z_{2,t^2q^4}\ar[d]^{2,t^3q^6}\ar[dr]^{1,t^2q^2}&&\alpha^L z_{2,t^2q^2}z_{1,t^3q^4}^{-1}z_{1,tq^2}\ar[dl]^{1,t^2q^4}\ar[d]^{2,t^3q^4}
\\\alpha^L z_{2,t^4q^8}^{-1}z_{1,t}\ar[d]^{1,t^2q^2} & z_{1,t^3q^6}^{-1}z_{1,t^3q^4}^{-1}z_{2,t^2q^4}z_{2,t^2q^2}
\ar[dl]^{2,t^3q^6}\ar[dr]^{2,t^3q^4}& \alpha^L z_{2,t^4q^6}^{-1}z_{1,tq^2}\ar[d]^{1,t^2q^4}
\\\alpha^L z_{1,t^3q^4}^{-1}z_{2,t^2q^2}z_{2,t^4q^8}^{-1}\ar[dr]^{2,t^3q^4}&&\alpha^L z_{1,t^3q^6}^{-1}z_{2,t^4q^6}^{-1}z_{2,t^2q^4}\ar[dl]^{2,t^3q^6}
\\&z_{2,t^4q^6}^{-1}z_{2,t^4q^8}^{-1}&}$$
 $\Pi_t^L$ gives the $t$-character of a tensor product of two $3$-dimensional fundamental
representation of $\U_t(A_3^{(1)})$ and $\Pi_q^L$ the $q$-character
of a $3$-dimensional fundamental representation of $\U_q(A_3^{(1)})$. 
$$\xymatrix{&Z_{1,-t}Z_{1,t}\ar[dl]^{1,t^2}\ar[dr]^{1,-t^2}&&Y_{1,q}\ar[dd]^{1,q^3}
\\ Z_{1,t^3}^{-1}Z_{1,-t}Z_{2,t^2}\ar[d]^{2,t^3}\ar[dr]^{1,-t^2}&& Z_{2,-t^2}Z_{1,-t^3}^{-1}Z_{1,t}\ar[dl]^{1,t^2}\ar[d]^{2,-t^3}&
\\  Z_{2,t^4}^{-1}Z_{1,-t}\ar[d]^{1,-t^2} & Z_{1,t^3}^{-1}Z_{1,-t^3}^{-1}Z_{2,t^2}Z_{2,-t^2}
\ar[dl]^{2,t^3}\ar[dr]^{2,-t^3}&  Z_{2,-t^4}^{-1}Z_{1,t}\ar[d]^{1,t^2}&Y_{1,q^5}^{-1}Y_{2,q^3}\ar[dd]^{2,q^5}
\\ Z_{1,-t^3}^{-1}Z_{2,-t^2}Z_{2,t^4}^{-1}\ar[dr]^{2,-t^3}&& Z_{1,t^3}^{-1}Z_{2,-t^4}^{-1}Z_{2,t^2}\ar[dl]^{2,t^3}&
\\&Z_{2,t^4}^{-1}Z_{2,-t^4}^{-1}&&Y_{2,q^7}^{-1}}$$
Consider $C_{i,j} = C_{j,i} = -1$ with $r_1^\vee = r_2^\vee = 2$. We choose $\phi(1) = 0$, $\phi(2) = 1$. We have
the following interpolating $(t,q)$-character.
$$\xymatrix{z_{2,1}\ar[d]^{2,qt}             & Z_{2,1}  \ar[d]^{2,t^2}               & Y_{2,1} \ar[d]^{2,q}
\\z_{2,q^2t^2}^{-1}z_{1,qt}\ar[d]^{1,q^2t^2} & Z_{2,t^4}^{-1}Z_{1,t^2}\ar[d]^{1,t^4} & Y_{2,q^2}^{-1}Y_{1,q} \ar[d]^{1,q^2}
\\z_{1,q^3t^3}^{-1}                          & Z_{1,t^6}^{-1}                        & Y_{1,q^3}^{-1}}$$
 $\Pi_t^L$ gives the $t$-character of a $3$-dimensional fundamental
representation of $\U_t(A_3^{(1)})$ and $\Pi_q^L$ the $q$-character
of a $3$-dimensional fundamental representation of $\U_q(A_3^{(1)})$. 

\subsection{Triple-laced case}

We use the notations of Section \ref{typeg}.
We need the function $\beta^L(t,q)$ such that $\beta^L(t,\epsilon) =
1$ and $\beta^L(1, q) = 0$ defined by $\beta^L(t,q) = 1 - \beta(q,t)$.

Consider the ring 
$$\Yim_{t,q}^L = \ZZ[X_{i,a}^{\pm 1}, \beta^L z_{i,a}^{\pm 1},
  \beta^L]_{i\in I, a\in\mathcal{C}}\subset \ZZ[z_{i,a}^{\pm
    1},\beta^L]_{i\in I, a\in\mathcal{C}},$$
\begin{equation*}
\begin{split}
\text{ where }X_{i,a} = \begin{cases} z_{i,a}&\text{ if $i = 2$,}
                     \\z_{i,aq^{-2}}z_{i,a}z_{i,aq^2}&\text{ if $i = 1$.}
\end{cases}
\end{split}
\end{equation*}

We then have the surjective specialization maps, $\Pi_t^L$, $\Pi_q^L$,
where for $a\in \mathcal{C}$, we use
$$X_{1,a} \mapsto Y_{1,a}\text{ , }z_{1,a} \mapsto Z_{1,-a}\text{ ,
}X_{2,a} \mapsto Y_{2,a}\text{ , }z_{2,a}\mapsto Z_{2,a^3}.$$ Note
that the identification is not one to one as $Z_{2,a^3}$ is identified
with $z_{2,a}$, $z_{2,\epsilon^2 a}$ and $Z_{2,\epsilon^4 a}$. Note
also that the identification is not involutive with respect to the
identification in Section \ref{typeg} as $z_{1,a}$ is identified with
$Z_{1,-a}$ and note with $Z_{1,a}$.

The ideal $\text{Ker}(\Pi_t^L)\cap \text{Ker}(\Pi_q^L)$ is generated
by the elements
$$\beta^L (\beta^L - 1)\text{ , }\beta^L (z_{i,aq} -
  z_{i,a\epsilon})\text{ , }\beta^L(z_{2,a} - z_{2,a\epsilon^2})
  \text{ , }(\beta^L - 1) (X_{i,at} - X_{i,a}),$$
  $$(X_{i,a} -
  X_{i,at})(X_{j,bq} - X_{j,b\epsilon})\text{ , }(X_{i,a} -
  X_{i,at})(z_{2,b} - z_{2,b\epsilon^2}),$$
for $i,j\in I$, $a,b\in\mathcal{C}$.

\begin{defi} For $a\in \mathcal{C}$ we define
$$\wt{B}_{1,a} = z_{i,a(tq^3)^{-1}}z_{i,atq^3} z_{2,a}^{-1}\text{ ,
  }\wt{B}_{2,a} = z_{2,a(qt)^{-1}}z_{2,aqt}
  z_{1,aq^{-2}}^{-1}z_{1,a}^{-1}z_{1,aq^2}^{-1}.$$
\end{defi}

\begin{lem} Let $a\in\mathcal{C}$. We have 

$\Pi_t^L(\wt{B}_{1,a}) = B_{1,(\Pi_t^L(a))}$, $\Pi_t^L(\wt{B}_{2,a}) =
  B_{2,-(\Pi_t^L(a))^3}$,

$\Pi_q^L(\wt{B}_{1,aq^{-2}}\wt{B}_{1,a}\wt{B}_{1,aq^2}) =
  A_{1,(\Pi_q^L(a))}$, $\Pi_q^L(\wt{B}_{2,a}) = A_{2,\Pi_q^L(a)}$.
\end{lem}

\demo Let $a' = \Pi_t^L(a)$. We have
$$\Pi_t^L(\wt{B}_{1,a}) = z_{1,-a' t^{-1}}z_{1,-a't} z_{2,a'}^{-1} =
Z_{1,a't^{-1}}Z_{1,a't}Z_{2,(a')^3}^{-1} = B_{1,a'},$$
$$\Pi_t^L(\wt{B}_{2,a}) = z_{2,-\epsilon^2 a' t^{-1}}z_{2,\epsilon a' t} z_{1,a'\epsilon^4}^{-1}z_{1,a'}^{-1}z_{1,a'\epsilon^2}^{-1}$$
$$= Z_{2,-(a')^3t^{-3}}Z_{2,-(a')^3t^3}
Z_{1,-a'\epsilon^4}^{-1}Z_{1,-a'}^{-1}Z_{1,-a'\epsilon^2}^{-1}
= B_{2,-(a')^3}.
$$
Let $a'' = \Pi_q^L(a)$. We have
$$\Pi_q^L((\wt{B}_{1,aq^{-2}}\wt{B}_{1,a}\wt{B}_{1,aq^2}))
=
(z_{1,a''q^{-5}}z_{1,a''q^{-3}}z_{1,a''q^{-1}})(z_{1,a''q}z_{1,a''q^3}z_{1,a''q^5})
z_{2,a''q^{-2}}^{-1}z_{2,a''}^{-1}z_{2,a''q^2}^{-1} $$
$$= Y_{1,a''q^{-3}}Y_{1,a''q^3} Y_{2,a''q^{-2}}^{-1}Y_{2,a''}^{-1}Y_{2,a''q^2}^{-1} = A_{1,a''},$$
$$\Pi_q^L(\wt{B}_{2,a}) = z_{2,a'' q^{-1}}z_{2,a''q}z_{1,a''}^{-1} = Y_{2,a''q^{-1}}Y_{2,a''}Y_{1,a''}^{-1} = A_{2,a''}.$$
\qed

Consider the subalgebra $\mathfrak{K}_{1,t,q}^L$ of $\Yim_{t,q}^L$
generated by the
$$X_{1,a}(1 + \beta^L\wt{B}_{1,q^5t}^{-1})
(1 + \beta^L\wt{B}_{1,q^3t}^{-1})
(1 + \beta^L\wt{B}_{1,qt}^{-1})
+(1 - \beta^L) X_{1,a}\wt{B}_{1,aq^5t}^{-1}
\wt{B}_{1,aq^3t}^{-1}\wt{B}_{1,aqt}^{-1},$$
$$\beta^L z_{1,a}(1 + \wt{B}_{1,aq^3t}^{-1}), z_{2,a}^{\pm 1}, \beta^L,$$ 
for $a\in\mathcal{C}$, and the subalgebra
$$\mathfrak{K}_{2,t,q}^L = \ZZ[z_{2,a}(1 + \wt{B}_{2,aqt}^{-1}),
  X_{1,a}^{\pm 1}, \beta^L z_{1,a}^{\pm 1}, \beta^L]_{a\in \mathcal{C}}.$$ 
Then we have the following.

\begin{lem} For $i\in I$, we have 
$\Pi_t^L(\mathfrak{K}_{i,t,q}^L)= \mathfrak{K}_{i,t}^L$ and 
$\Pi_q^L(\mathfrak{K}_{i,t,q}^L) = \mathfrak{K}_{i,q}$.
\end{lem}

\demo $\Pi_t^L(\mathfrak{K}_{1, t,q}^L)$ is equal to
$$\ZZ[z_{1,\epsilon^{-2} a}z_{1,a}z_{1,a\epsilon^2}(1 + 
B_{1,-\epsilon^2 at}^{-1})(1 + 
B_{1,-at}^{-1})(1 + 
B_{1,-\epsilon^{-2} at}^{-1}),
z_{1,a}(1 + B_{1,- a t}^{-1}), z_{2,a}^{\pm 1}]_{a\in t^{\ZZ}\epsilon^\ZZ}$$
$$= \ZZ[Z_{1,-a}(1 + B_{1,-a t}^{-1}), Z_{2,a^3}^{\pm 1}]_{a\in
  t^{\ZZ}\epsilon^\ZZ} = \mathfrak{K}_{1,t}^L.$$
We also have
$$\Pi_q^L(\mathfrak{K}_{1, t,q}^L) = \ZZ[\Pi_q^L(X_{1,a}(1 +
  \wt{B}_{i,aq^5t}^{-1}\wt{B}_{i,aq^3t}^{-1}\wt{B}_{i,aqt}^{-1})),
\Pi_q^L(X_{2,a}^{\pm 1})]_{a\in\mathcal{C}}$$
$$= \ZZ[Y_{1,a}(1 + A_{1,aq^3}^{-1}), Y_{2,a}^{\pm 1}]_{a\in q^{\ZZ}}
= \mathfrak{K}_{1,q},$$
$$\Pi_t^L(\mathfrak{K}_{2,t,q}^L)
=\ZZ[z_{2,a}(1 + B_{2,a^3 t^3}^{-1}),
  z_{1,a}^{\pm 1}]_{a\in t^\ZZ\epsilon^\ZZ}
  = \ZZ[Z_{2,a^3}(1 + B_{2,a^3 t^3}^{-1}),
  Z_{1,a}^{\pm 1}]_{a\in t^\ZZ\epsilon^\ZZ} = \mathfrak{K}_{2,t}^L,$$ 
and we have
$$\Pi_q^L(\mathfrak{K}_{2,t,q}^L) = \ZZ[Y_{2,a}(1 + A_{2,aq}^{-1}),
  Y_{1,a}^{\pm 1}]_{a\in q^\ZZ} = \mathfrak{K}_{2,q}.$$ 
\qed

We set
\begin{equation*}
(\Yim_t^L)' = \ZZ[Z_{2,a}^{\pm 1}]_{a\in (\epsilon^\ZZ t^\ZZ)^3}\otimes
\ZZ[(Z_{1,a}Z_{1,a\epsilon^2}Z_{1,a\epsilon^4})^{\pm 1}]_{a\in
  \epsilon^\ZZ t^\ZZ}\subset \Yim_t^L,
\end{equation*}
and we define $\mathfrak{K}^L\subset \tilde{\Yim}_{t,q}^L$ as above.

As in the previous sections, we check the existence of elements in
$\mathfrak{K}^L$ that we call interpolating $(t,q)$-characters. First,
we have the following.

\begin{figure}[htbp]
\begin{center}
\leavevmode
\begin{equation*}
\divide\dgARROWLENGTH by 3
\dgARROWPARTS=6
\begin{diagram}
  \node[3]{2_1} \arrow{s,r}{2,qt}
\\
  \node[3]{2_{q^2t^2}^{-1} 1_{q^{-1} t} 1_{qt} 1_{q^3 t}}
  \arrow{sw,t,3}{1,q^2 t^2} \arrow{s,l,3}{1,q^4 t^2}
  \arrow{se,t,3}{1,q^6 t^2}
\\
  \node[2]{\beta^L 1_{q^5 t^3}^{-1} 1_{qt} 1_{q^3 t}}
  \arrow{s,l}{1,q^4 t^2} \arrow{se,b,1}{1,q^6t^2}
  \node{\beta^L 1_{q^{-1} t} 1_{q^7t^3}^{-1} 1_{q^3 t}}
  \arrow{sw,b,1}{1,q^2 t^2}\arrow{se,b,1}{1, q^6 t^2}
  \node{\beta^L 1_{q^{-1} t}1_{q t} 1_{q^9 t^3}^{-1}}
  \arrow{sw,b,1}{1, q^2 t^2} \arrow{s,r}{1, q^4 t^2}
\\
  \node[2]{\beta^L 1_{q^5 t^3}^{-1} 2_{q^4t^2} 1_{q^7t^3}^{-1} 1_{q^3 t}}
  \arrow{sw,b,1}{2,q^5t^3} \arrow{se,b,1}{1,q^6 t^2}
  \node{\beta^L 1_{q^5 t^3}^{-1} 2_{q^6 t^2} 1_{q t} 1_{q^9 t^3}^{-1}}
  \arrow{sw,b,1}{2,q^7t^3}
  \arrow{s,r}{1,q^4 t^2}
  \node{\beta^L 1_{q^{-1} t} 2_{q^6 t^2} 1_{q^7 t^3}^{-1} 1_{q^9 t^3}^{-1}}
  \arrow{sw,b,1}{1, q^2 t^2}
  \arrow{s,r}{2, q^7t^3}
\\
  \node{\beta^L 2_{q^6 t^4}^{-1}1_{q^3 t}1_{q^3 t^3}}
  \arrow{s,l}{1,q^6 t^4}
  \node{\beta^L 1_{qt} 1_{q t^3} 2_{q^2t^4}^{-1}}
  \arrow{s,r,1}{1, q^{10} t^4}
  \node{2_{q^4t^2} 2_{q^6t^2} 1_{q^5 t^3}^{-1} 1_{q^7 t^3}^{-1} 1_{q^9 t^3}^{-1}}
  \arrow{s,r}{2, q^5 t^3}
  \node{\beta^L 2_{q^8t^4}^{-1} 1_{q^{-1} t} 1_{q^{-1} t^3}}
  \arrow{s,r}{1, q^8 t^4}
\\
  \node{\beta^L 1_{q^3 t}1_{q^9 t^5}^{-1}}
  \arrow{s,l}{1, q^6 t^2}
  \node{\beta^L 1_{qt}1_{q^7 t^5}^{-1}}
  \arrow{s,l}{1, q^4 t^2}
  \node{(1 + \beta^L) 2_{q^{10}t^4}^{-1} 2_{q^6t^2}}
  \arrow{s,r}{2, q^7 t^3}
  \node{\beta^L 1_{q^{-1} t} 1_{q^5 t^5}^{-1}}\arrow{s,r}{1, q^2 t^2}
\\
  \node{\beta^L 1_{q^9 t^3}^{-1}1_{q^9 t^5}^{-1}2_{q^6t^2}}
  \arrow{se,b,1}{2, q^7t^3}
  \node{\beta^L 1_{q^7 t^3}^{-1}1_{q^7 t^5}^{-1} 2_{q^4 t^2}}
  \arrow{se,b,1}{2, q^5 t^3}
  \node{1_{q^3 t^3} 2_{q^6t^4}^{-1}2_{q^8t^4}^{-1} 1_{q^5 t^3}1_{q^7 t^3}}
  \arrow{sw,b,1}{1, q^6 t^4}\arrow{s,r,2}{1, q^{10} t^4}
  \arrow{se,t,3}{1, q^8 t^4}
  \node{\beta^L 2_{q^4t^2}1_{q^5 t^3}^{-1}1_{q^5 t^5}^{-1}}
  \arrow{s,r}{2, q^5 t^3}
\\
  \node[2]{\beta^L 1_{q^5 t^3}  2_{q^8t^4}^{-1}1_{q^9 t^5}^{-1}1_{q^7 t^3}
  }
  \arrow{s,l}{1, q^{10} t^4}
  \arrow{se,b,1}{1,q^8 t^4}
  \node{\beta^L 1_{q^5 t^3}  2_{q^8 t^4}^{-1} 1_{q^3 t^3} 1_{q^{13} t^5}^{-1}
  }
  \arrow{sw,b,1}{1,q^6 t^4}\arrow{se,b,1}{1, q^8 t^4}
  \node{\beta^L 1_{q^{11}t^5}^{-1} 2_{q^{10}t^4}^{-1} 1_{q^3 t^3} 1_{q^7 t^3}}
  \arrow{sw,b,1}{1, q^6 t^4} \arrow{s,r}{1, q^{10} t^4}
\\
  \node[2]{\beta^L 1_{q^5 t^3} 1_{q^9 t^5}^{-1} 1_{q^{13} t^5}^{-1}}
  \arrow{se,b,2}{1, q^8 t^4}
  \node{\beta^L 1_{q^9 t^5}^{-1} 1_{q^{11} t^5}^{-1} 1_{q^7 t^3}}
  \arrow{s,r}{1, q^{10} t^4}
  \node{\beta^L 1_{q^{11} t^5}^{-1} 1_{q^3 t^3} 1_{q^{13} t^5}^{-1}}
  \arrow{sw,b,2}{1, q^6 t^4}
\\
  \node[3]{1_{q^9 t^5}^{-1}2_{q^{10}t^4}1_{q^{11} t^5}^{-1}1_{q^{13} t^5}^{-1}}
  \arrow{s,r}{2, q^{11}t^5}
\\
  \node[3]{2_{q^{12}t^6}^{-1}}
\end{diagram}
\end{equation*}
\end{center}
\end{figure}

\newpage

\begin{rem} In the diagram we use the following identities which hold in
$\tilde{\Yim}_{t,q}$:
$$\beta^L 2_{a} = \beta^L 2_{aq^2}\text{ and }\beta^L 1_a = \beta^L
1_{aq^6}.$$
\end{rem}

It is easy to check that the expression is in the
$\mathfrak{K}^L$. For example,
$$2_{q^4t^2} 2_{q^6t^2} 1_{q^5 t^3}^{-1} 1_{q^7 t^3}^{-1} 1_{q^9 t^3}^{-1}
+ (1 + \beta^L) 2_{q^{10}t^4}^{-1} 2_{q^6t^2}
+ 1_{q^3 t^3} 2_{q^6t^4}^{-1}2_{q^8t^4}^{-1} 1_{q^5 t^3}1_{q^7 t^3}
$$
$$= (\beta^L - 1)  + 2_{q^4t^2} 2_{q^6t^2} 1_{q^5 t^3}^{-1} 1_{q^7 t^3}^{-1} 1_{q^9 t^3}^{-1} (1 + \wt{B}_{2,q^5t^3})(1 + \wt{B}_{2,q^7t^3})\in \mathfrak{K}_{2,t,q}^L.$$
The image by $\Pi_t^L$ is the twisted $t$-characters of a
29-dimensional fundamental representation of $\U_q(D_4^{(3)})$ (see
the diagram in Section \ref{trex}).

\noindent The image by $\Pi_q^L$ is the following $q$-character of
a $7$-dimensional fundamental representation of $\U_q(G_2^{(1)})$, as computed in \cite[Appendix]{hadv}.
$$2_1 + 2_{q^2}^{-1}1_q + 1_{q^7}^{-1}2_{q^4}2_{q^6} +
2_{q^4}2_{q^8}^{-1} + 2_{q^6}^{-1}2_{q^8}^{-1}1_{q^5} +
1_{q^{11}}^{-1}2_{q^{10}} + 2_{q^{12}}^{-1}.$$

Now we have to consider the case of the monomial
$Z_{1,\epsilon}Z_{1,\epsilon^3}Z_{1,\epsilon^5}$. The dimension of the
corresponding simple module of $\U_t(D_4^{(3)})$ is
$8^3 = 512$ (the module is the tensor product of three fundamental
representations of dimension $8$).

There is also an interpolating $(t,q)$-character in this case. We do
not list all $512$ monomials, but we list the $15$ monomials which do 
not have $\beta^L$ in their coefficient:

$1_1 1_{q^2}1_{q^4}; 
1_{q^6t^2}^{-1} 1_{q^8 t^2}^{-1} 1_{q^{10}t^2}^{-1} 2_{q^3t} 2_{q^5t} 2_{q^7t};
2_{q^3t} 2_{q^5t}2_{q^9t^3}^{-1};
1_{q^4t^2}1_{q^6t^2}1_{q^8t^2}2_{q^3t}2_{q^7t^3}^{-1}2_{q^9t^3}^{-1};
\\1_{q^2t^2}1_{q^4t^2}^21_{q^6t^2}^21_{q^8t^2}2_{q^5t^3}^{-1}2_{q^7t^3}^{-1}2_{q^9t^3}^{-1};
1_{q^{10}t^4}^{-1}1_{q^{12}t^4}^{-1}1_{q^{14}t^4}^{-1}2_{q^3t}2_{q^{11}t^3};
\\1_{q^2t^2}1_{q^4t^2}1_{q^6t^2}1_{q^{10}t^4}^{-1}1_{q^{12}t^4}^{-1}1_{q^{14}t^4}^{-1}2_{q^5t^3}^{-1}2_{q^{11}t^3};
1_{q^8t^4}^{-1}1_{q^{10}t^4}^{-1}1_{q^{12}t^4}^{-1}1_{q^4t^2}1_{q^6t^2}1_{q^8t^2};
2_{q^3t}2_{q^{13}t^5}^{-1};
\\1_{q^8t^4}^{-1}1_{q^{10}t^4}^{-2}1_{q^{12}t^4}^{-2}1_{q^{14}t^4}^{-1}2_{q^7t^3}2_{q^9t^3}2_{q^{11}t^3};
2_{q^5t^3}^{-1}2_{q^{13}t^5}^{-1}1_{q^2t^2}1_{q^4t^2}1_{q^6t^2};
2_{q^7t^3}2_{q^9t^3}2_{q^{13}t^5}^{-1}1_{q^8t^4}^{-1}1_{q^{10}t^4}^{-1}1_{q^{12}t^4}^{-1};
\\2_{q^7t^3}2_{q^{11}t^5}^{-1}2_{q^{13}t^5}^{-1};
2_{q^9t^5}^{-1}2_{q^{11}t^5}^{-1}2_{q^{13}t^5}^{-1}1_{q^6t^4}1_{q^8t^4}1_{q^{10}t^4};
1_{q^{12}t^6}^{-1}1_{q^{14}t^6}^{-1}1_{q^{16}t^6}^{-1}.
$

As the other terms disappear when we specialize at $t = 1$, we
can compute the specialization from the above terms. We get
the $q$-character of a $15$-dimensional fundamental
representation of $\U_q(G_2^{(1)})$ (image by $\Pi_q$ of the first
example in Section \ref{trex}).

\subsection{Conclusion}
We go back to the general case, that is, $r=2$ or $r =3$.

With the existence of the elementary interpolating
$(t,q)$-characters in the subsections, we can conclude the proof of
the main results of this section.

\begin{thm}\label{virtrev} For all dominant $m$ there is a unique
  $F(m)\in\mathfrak{K}^L$ such that $m$ is the unique dominant
  monomial of $F(m)$.\end{thm}

We state its consequence in terms of KR modules. 

As for $i\in I_1^\vee$, the only KR module of
$\U_t({}^L\hat{\Glie})$ for the node $i$ with highest monomial in
$(\Yim_t^L)'$ is trivial, we extend the definition. For $i\in
I_1^\vee$, a simple $\U_t({}^L\wh{\Glie})$-module with the highest
monomial of the form
$$(Z_{i,a}Z_{i,at^2}\cdots
Z_{i,at^{2(k-1)}})(Z_{i,-a}Z_{i,-at^2}\cdots Z_{i,at^{2(k-1)}})$$
for the double-laced case and of the form
$$(Z_{1,a}Z_{1,at^2}\cdots Z_{1,at^{2(k-1)}})(Z_{1,\epsilon^2
  a}Z_{1,\epsilon^2 at^2}\cdots Z_{1,\epsilon^2
  at^{2(k-1)}})(Z_{1,\epsilon^4 a}Z_{1,\epsilon^4 at^2}\cdots
Z_{1,\epsilon^4 at^{2(k-1)}})$$
for the triple-laced case, will also be called a KR module.

As above, we have the notion of Langlands dual
representation and interpolating $(t,q)$-character $\chi\in \mathfrak{K}^L$
with highest monomial in $\wt{\Yim}_{t,q}\setminus \alpha^L
\wt{\Yim}_{t,q}$ ($\wt{\Yim}_{t,q}\setminus \beta^L \wt{\Yim}_{t,q}$
in the triple-laced case): $\Pi_t^L(\chi)$ is the twisted $t$-character of
a $\U_t({}^L\hat{\Glie})$-module and $\Pi_q^L(\chi)$ is 
the $q$-character of a $\U_q(\hat{\Glie})$-module. 
We obtain the following consequence of Theorem \ref{virtrev}:

\begin{thm}\label{dualrev} An irreducible tensor product of KR modules
  over $\U_q({}^L\hat{\Glie})$
  of highest monomial $M\in(\Yim_t^L)'$ has a Langlands dual
  representation. Moreover, the Langlands dual representation of a KR
  module over $\U_q({}^L\hat{\Glie})$ is a KR module over
  $\U_q(\hat{\Glie})$.
\end{thm}


\begin{thebibliography}{99}

\bibitem{bpil}{\bf P. Bouwknegt and K. Pilch}, {\it On deformed
  $W$-algebras and quantum affine algebras},
  {Adv. Theor. Math. Phys. {\bf 2} (1998),  no. 2, 357--397}.

\mk

\bibitem{ch}{\bf V. Chari and D. Hernandez}, {\it Beyond
  Kirillov--Reshetikhin modules}, {Contemp. Math. {\bf 506} (2010), 49--81}.

\mk

\bibitem{cp}{\bf V. Chari and A. Pressley}, {\it A Guide to Quantum
  Groups}, {Cambridge University Press, Cambridge, 1994}.

\mk

\bibitem{cp2}{\bf V. Chari and A. Pressley}, {\it Factorization of
    representations of quantum affine algebras}, {Modular interfaces,
    (Riverside CA 1995), AMS/IP Stud. Adv. Math., 4 (1997), 33-40}

\mk

\bibitem{FH}{\bf E. Frenkel and D. Hernandez}, {\it Langlands duality
  for representations of quantum groups}, {to appear in Math. Ann.}

\mk

\bibitem{fm}{\bf E. Frenkel and E. Mukhin}, {\it Combinatorics of
$q$-characters of finite-dimensional representations of quantum affine
algebras}, {Comm. Math. Phys. {\bf 216} (2001), no. 1, 23--57}.

\mk

\bibitem{fr2}{\bf E. Frenkel and N. Reshetikhin}, {\it Deformations of
$W$-algebras associated to simple Lie algebras},
{Comm. Math. Phys. {\bf 197} (1998), no. 1, 1--32}.

\mk

\bibitem{fr3}{\bf E. Frenkel and N. Reshetikhin}, {\it The
$q$-characters of representations of quantum affine algebras and
deformations of $W$-algebras}, {Recent developments in quantum affine
algebras and related topics (Raleigh, NC, 1998), pp. 163--205,
Contemp. Math., 248, Amer. Math. Soc., Providence, RI, 1999}.

\mk

\bibitem{hadv}{\bf D. Hernandez}, {\it Algebraic Approach to
  q,t-Characters}, {Adv. Math. {\bf 187} (2004), no. 1, 1--52}.

\mk

\bibitem{haja}{\bf D. Hernandez}, {\it The $t$-Analogs of
$q$-Characters at Roots of Unity for Quantum Affine Algebras and
Beyond}, {J. Algebra {\bf 279} (2004), no. 2 , 514--557}.

\mk

\bibitem{hcr}{\bf D. Hernandez}, {\it The Kirillov--Reshetikhin
conjecture and solutions of $T$-systems}, {J. Reine Angew. Math. {\bf
596} (2006), 63--87}.

\mk

\bibitem{hcmp}{\bf D. Hernandez}, {\it On minimal affinizations of
  representations of quantum groups},  {Comm. Math. Phys. {\bf 277}
  (2007), 221--259.}

\mk

\bibitem{hgen}{\bf D. Hernandez}, {\it Kirillov--Reshetikhin
  conjecture : the general case}, 
  {Int. Math. Res. Not. {\bf 2010}, no. 1, 149--193}.

\mk

\bibitem{hprod}{\bf D. Hernandez}, {\it Simple tensor products},
  {Invent. Math. {\bf 181} (2010), no. 3, 649--675}.

\mk

\bibitem{hl}{\bf D. Hernandez and B. Leclerc}, {\it Cluster algebras and quantum affine algebras}, 
{Duke Math. J. {\bf 154} (2010), no. 2, 265--341}.

\mk

\bibitem{knh} {\bf A. Kuniba, S. Nakamura and R. Hirota}, {\it
Pfaffian and determinant solutions to a discretized Toda equation for
$B_r$, $C_r$ and $D_r$}, {J. Phys. A {\bf 29} (1996), no. 8, 1759--1766}.

\mk

\bibitem{kos} {\bf A. Kuniba, Y. Ohta and J. Suzuki}, {\it Quantum
Jacobi-Trudi and Giambelli Formulae for $\mathcal{U}_q(B_r^{(1)})$
from Analytic Bethe Ansatz}, {J. Phys. A {\bf 28} (1995), no. 21, 6211--6226}.

\mk

\bibitem{ks}{\bf A. Kuniba and S. Suzuki}, {\it Analytic Bethe Ansatz
for fundamental representations of Yangians},
{Commun. Math. Phys. {\bf 173} (1995), 225 - 264}.

\mk

\bibitem{pl}{\bf P. Littelmann}, {\it The path model, the quantum Frobenius map and standard monomial theory}, 
{Algebraic groups and their representations (Cambridge, 1997), 175--212,
NATO Adv. Sci. Inst. Ser. C Math. Phys. Sci., 517, Kluwer Acad. Publ., Dordrecht, 1998}

\mk

\bibitem{M} {\bf K. McGerty}, {\em Langlands duality for
  representations and quantum groups at a root of unity}, 
  {Comm. Math. Phys. {\bf 296}  (2010),  no. 1, 89--109}.

\mk

\bibitem{Nab} {\bf H. Nakajima}, {\it Quiver Varieties and $t$-Analogs
of $q$-Characters of Quantum Affine Algebras}, {Ann. of Math. {\bf
160} (2004), 1057 - 1097}.

\mk

\bibitem{Nad}{\bf H. Nakajima}, {\it $t$-analogs of $q$-characters of
Kirillov--Reshetikhin modules of quantum affine algebras},
{Represent. Theory {\bf 7} (2003), 259--274 (electronic)}.

\end{thebibliography}
\end{document}